\documentclass{article}

\usepackage{times}
\usepackage{graphicx} 
\usepackage{subfigure} 

\usepackage{natbib}

\usepackage{hyperref}

\usepackage{amsfonts}
\usepackage{mathrsfs}
\usepackage{amsmath, amsthm}
\usepackage{amssymb,url,paralist}

\newtheorem{thm}{Theorem}

\newtheorem{lemma}{Lemma}

\usepackage[accepted]{icml2015}

\def \R {\mathbb{R}}

\def \E {\mathrm{E}}
\def \x {\mathbf{x}}

\def \v {\mathbf{v}}

\def \X {\mathcal{X}}

\def \xh {\widehat{\x}}

\def \u {\mathbf{u}}

\def \v {\mathbf{v}}

\def \R {\mathbb{R}}

\def \b {\mathbf{b}}

\def \s {\mathbf{s}}

\icmltitlerunning{Sampling Dependent Spectral Error Bound for CSS}

\begin{document} 

\twocolumn[
\icmltitle{An Explicit Sampling Dependent Spectral Error Bound\\ for Column Subset Selection}

\icmlauthor{Tianbao Yang}{tianbao-yang@uiowa.edu}
\icmladdress{Department of Computer Science, the University of Iowa, Iowa City, USA}
\icmlauthor{Lijun Zhang}{zhanglj@lamda.nju.edu.cn}
\icmladdress{National Key Laboratory for Novel Software Technology, Nanjing University, Nanjing, China}
\icmlauthor{Rong Jin}{rongjin@cse.msu.edu}
\icmladdress{Department of Computer Science and Engineering, Michigan State University, East Lansing, USA\\Institute of Data Science and Technologies at Alibaba Group, Seattle, USA}
\icmlauthor{Shenghuo Zhu}{shenghuo@gmail.com}
\icmladdress{Institute of Data Science and Technologies at Alibaba Group, Seattle, USA}
\icmlkeywords{boring formatting information, machine learning, ICML}

\vskip 0.3in
]

\begin{abstract} 
In this paper, we consider the problem of column subset selection.  We present a novel analysis of the spectral norm reconstruction for a simple randomized algorithm and establish a new bound that depends explicitly on the sampling probabilities. The sampling dependent error bound (i) allows us to better understand the tradeoff in the reconstruction error  due to sampling probabilities, (ii) exhibits more insights than existing error bounds that exploit  specific probability distributions, and (iii) implies better sampling distributions. In particular, we show that a sampling distribution with probabilities  proportional to the square root of the statistical leverage scores is always better than uniform sampling and is better than leverage-based sampling when the statistical leverage scores are very nonuniform. And  by solving a constrained optimization problem related to the error bound with an efficient bisection search  we are able to achieve better performance than using either the leverage-based distribution or that proportional to  the square root of the statistical leverage scores. Numerical simulations  demonstrate the benefits of the new sampling distributions for low-rank matrix approximation and least square approximation compared to state-of-the art algorithms. 
\end{abstract}
\setlength{\belowdisplayskip}{2pt} \setlength{\belowdisplayshortskip}{2pt}
\setlength{\abovedisplayskip}{2pt} \setlength{\abovedisplayshortskip}{2pt}

\section{Introduction}
Give a  data matrix $A\in\R^{m\times n}$, {\bf column subset selection (CSS)} is an important technique for constructing a compressed representation and a low rank approximation of $A$ by selecting  a small number of columns. Compared with conventional singular value decomposition (SVD), CSS could  yield more interpretable output while maintaining performance as close as SVD~\cite{mahoney-2011-randomized}. Recently, CSS has been applied successfully to problems of interest to geneticists such as genotype reconstruction, identifying substructure in heterogeneous populations, etc.~\cite{citeulike:1682474,pmid17151345,pmid20805874,pmid21902678}.  

Let $C\in\R^{m\times\ell}$ be the matrix formed by $\ell$ selected columns of $A$. The key question to CSS is how to select the columns to minimize the reconstruction error: 
\begin{align*}
\|A - P_CA\|_\xi,
\end{align*}
where $P_C=CC^{\dagger}$ denotes the projection onto the column space of $C$ with $C^{\dagger}$ being the pseudo inverse of $C$ and $\xi=2$ or $F$ denotes the spectral norm or the Frobenius norm. In this paper, we are particularly interested in the spectral norm reconstruction  with respect to a target rank $k$. 

Our analysis is based on a randomized algorithm that selects $\ell>k$ columns from $A$ according to  sampling probabilities $\s = (s_1, \ldots, s_n)$. Building on advanced matrix concentration inequalities (e.g., matrix Chernoff bound and Bernstein inequality), we develop a novel analysis of the spectral norm reconstruction and establish  a sampling dependent relative spectral error bound with a high probability as following: 
\begin{align*}
\|A - P_CA\|_2\leq (1+\epsilon(\s))\|A - A_k\|_2,
\end{align*}
where $A_k$ is the best rank-$k$ approximation of $A$ based on SVD and $\epsilon(\s)$ is a quantity dependent on the sampling probabilities $\s$ besides the scalars $n, k, \ell$.  As revealed in our main theorem (Theorem~\ref{thm:final-bound}), the quantity  $\epsilon(\s)$ also depends on the {\bf statistical leverage scores (SLS)} inherent to the data, based on which are several important randomized algorithms for CSS. 

To the best of our knowledge, this is the first such kind of error bound for CSS. Compared with existing error bounds, the sampling dependent error bound brings us several benefits: (i) it allows us to better understand the tradeoff in the spectral error of reconstruction due to sampling probabilities, complementary to a recent result on the tradeoff  from a statistical perspective~\cite{DBLP:conf/icml/MaMY14} for least square regression; (ii) it implies that a distribution with sampling probabilities proportional to the square root of the SLS is always better than the uniform sampling, and is potentially better than that proportional to the SLS when they are skewed; (iii) it motivates an optimization approach by solving a constrained optimization problem related to the error bound to attain better performance. In addition to the theoretical analysis,  we also develop an efficient bisection search algorithm to solve the constrained optimization problem for finding better sampling probabilities.  

By combining  our analysis with recent developments for spectral norm reconstruction of CSS~\cite{conf/focs/BoutsidisDM11}, we also establish the same error bound for an exact rank-$k$ approximation, i.e., 
\begin{align*}
\left\|A - \Pi^2_{C,k}(A)\right\|_2\leq (1+\epsilon(\s))\|A - A_k\|_2,
\end{align*}
where $\Pi^2_{C,k}(A)$ is the best approximation to $A$ within the column space of $C$ that has rank at most $k$. 

The remainder of the paper is organized as follows. We review some closely related work in~Section~\ref{sec:related}, and present the main result in Section~\ref{sec:main} with some preliminaries in Section~\ref{sec:preliminary}. We  conduct  some empirical studies in Section~\ref{sec:empirical} and present the detailed analysis  in Section~\ref{sec:analysis}. Finally, conclusion is made. 

\section{Related Work}\label{sec:related}
In this section, we review some previous  work on CSS, low-rank matrix approximation, and other closely related work on randomized algorithms for matrices. We focus our discussion on the spectral norm reconstruction. 

Depending on whether the columns are selected deterministically or randomly, the algorithms for CSS can be categorized into deterministic algorithms and randomized algorithms. Deterministic algorithms select $\ell\geq k$ columns with some deterministic selection criteria. Representative algorithms in this category are rank revealing QR factorization and its variants  from the filed of numerical linear algebra~\cite{Gu:1996:EAC:239729.239737,Pan2000199,pan99}. A recent work~\cite{conf/focs/BoutsidisDM11} based on the dual set spectral sparsification also falls into this category which will be discussed shortly.  Randomized algorithms usually define  sampling probabilities $\s\in\R^n$ and then select $\ell\geq k$ columns based on these sampling probabilities. Representative sampling probabilities include ones that depend  the squared Euclidean norm of columns (better for Frobenius norm reconstruction)~\cite{Frieze:2004:FMA:1039488.1039494}, the squared volume of simplices defined by the selected subsets of columns (known as volume sampling)~\cite{journals/corr/abs-1004-4057}, and the SLS (known as {\bf leverage-based sampling} or subspace sampling)~\cite{journals/siammax/DrineasMM08,Boutsidis:2009:IAA:1496770.1496875}. 

Depending on whether $\ell>k$ is allowed, the error bounds for CSS  are different. Below, we review several representative error bounds.   If exactly $k$ columns are selected to form $C$, the best bound was achieved by the rank revealing QR factorization~\cite{Gu:1996:EAC:239729.239737} with the error bound given by:
\begin{align}
\|A - P_CA\|_2\leq \sqrt{1 + k(n-k)}\|A - A_k\|_2.
\end{align}
with a running time $O(mnk\log n)$. The same error  bound was also achieved by using volume sampling~\cite{journals/corr/abs-1004-4057}. The running time of volume sampling based algorithms can be made close to linear to the size of the target matrix.  \citet{Boutsidis:2009:IAA:1496770.1496875} proposed a two-stage algorithm for selecting exactly $k$ columns and provided error bounds for both the spectral norm and the Frobenius norm, where in the firs stage $\Theta(k\log k)$ columns are  sampled based on a distribution related to the SLS and more if for the spectral norm reconstruction and in the second stage $k$ columns are selected based on the rank revealing QR factorization. The spectral error bound in this work that holds with a constant probability 0.8 is following: 
\begin{align}\label{eqn:bound-0}
\|A - &P_CA\|_2\leq\\
& \Theta \left(k\log^{1/2}k +n^{1/2}k^{3/4}\log^{1/4}(k)  \right)\|A - A_k\|_2\nonumber
\end{align}
The time complexity of their algorithm (for the spectral norm reconstruction) is given by $O(\min(mn^2, m^2n))$ since it requires SVD of the target matrix for computing the sampling probabilities. 

If more than $k$ columns are allowed to be selected, i.e., $\ell>k$, better error bounds can be achieved. In the most recent work by~\citet{conf/focs/BoutsidisDM11}, nearly optimal error bounds were shown by selecting $\ell>k$ columns with a deterministic selection criterion based on the dual set spectral sparsification.  In particular, a deterministic polynomial-time algorithm~\footnote{A slower deterministic algorithm with a time complexity $T_{\text{SVD}} + O(\ell n(k^2 + (\rho -k)^2))$ was also presented with an error bound $O(\sqrt{\rho/\ell})\|A - A_k\|_2$, where $\rho$ is the rank of $A$.} was proposed that achieves the following error bound:
\begin{align}\label{eqn:bound-2}
\|A - P_{C}A\|_2&\leq \left(1 + \frac{1+\sqrt{n/\ell}}{1 - \sqrt{k/\ell}}\right)\|A- A_k\|_2
\end{align}
 in $T_{V_k} + O(n\ell k^2)$ time where $T_{V_k}$ is the time needed to compute the top $k$ right singular vectors of $A$ and $O(n\ell k^2)$ is the time needed to compute the selection scores. This bound is close to the lower bound $\Omega\left(\sqrt{\frac{n+\alpha^2}{\ell+\alpha^2}}\right), \alpha>0$ established in their work. It is worth mentioning  that the selection scores in~\cite{conf/focs/BoutsidisDM11} computed based on the dual set spectral sparsification  is difficult to understand than the SLS. 
 
 Although our sampling dependent error bound is not directly comparable to these results, our analysis exhibits that the derived error bound could be much better than that in~(\ref{eqn:bound-0}). When the SLS are nonuniform, our new sampling distributions could lead to a better result than~(\ref{eqn:bound-2}).  Most importantly, the sampling probabilities in our algorithm are only related to the SLS and that can be  computed more efficiently (e.g.,  exactly in $O(T_{V_k})$ or approximately in $O(mn\log n)$~\cite{DBLP:journals/jmlr/DrineasMMW12}). In simulations, we observe that the new sampling distributions could yield even better  spectral norm reconstruction  than the deterministic selection criterion in~\cite{conf/focs/BoutsidisDM11}, especially when the SLS are nonuniform.

For low rank matrix approximation, several other randomized algorithms have been  recently developed. For example, \citet{Halko:2011:FSR:2078879.2078881} used a random Gaussian matrix $\Omega\in\R^{n\times \ell}$ or a subsampled random Fourier transform to construct a matrix $\Omega$ and then project $A$ into the column space of $Y=A\Omega$, and they established numerous spectral error bounds. Among them is a comparable error bound $O(\sqrt{n/\ell})\|A - A_k\|_2$ to~(\ref{eqn:bound-2}) using the subsampled random Fourier transform. Other randomized algorithm for low rank approximation include CUR decomposition~\cite{journals/siammax/DrineasMM08,NIPS2012_0301,Wang:2013:ICM:2567709.2567748} and the Nystr\"{o}m based approximation for PSD matrices~\cite{Drineas05onthe,conf/icml/GittensM13}. 

Besides  low rank matrix approximation and column selection, CSS has also been successfully applied to least square approximation, leading to faster and interpretable  algorithms for over-constrained least square regression. In particular, if let $\Omega\in\R^{\ell\times m}$ denote a scaled sampling matrix  corresponding to selecting $\ell<m$ rows from $A$, the least square problem $\min_{\x \in\R^n}\|A\x - \b\|^2_2$ can be approximately solved by $\min_{\x\in\R^n}\|\Omega A\x - \Omega\b\|^2_2$~\cite{journals/siammax/DrineasMM08,conf/soda/DrineasMM06,Drineas:2011:FLS:1936922.1936925}. At \textit{ICML 2014}, \citet{DBLP:conf/icml/MaMY14} studied CSS for least square approximation from a statistical perspective. They showed the expectation and  variance of the solution to the approximated least square with uniform sampling and leverage-based sampling. They found that leveraging based estimator could suffer from a large variance when the SLS are very nonuniform while uniform sampling is less vulnerable to very small SLS. This tradeoff is complementary to our observation. However, our observation follows directly from the spectral norm error bound. Moreover, our analysis reveals  that the sampling distribution with probabilities proportional to the square root of the SLS is always better than uniform sampling, suggesting  that intermediate sampling probabilities between SLS and their square roots by solving a constrained  optimization problem could yield better performance than the mixing strategy that linearly combines the SLS and uniform probabilities as suggested in~\cite{DBLP:conf/icml/MaMY14}.

There are much more work on studying the Frobenius norm reconstruction of CSS~\cite{conf/approx/DrineasMM06,Guruswami:2012:OCL:2095116.2095211,conf/focs/BoutsidisDM11,journals/siammax/DrineasMM08,Boutsidis:2009:IAA:1496770.1496875}. For more references, we refer the reader to the survey~\cite{mahoney-2011-randomized}.  It remains an interesting question to establish  sampling dependent error bounds for other randomized matrix algorithms.

\vspace*{-0.1in}
\section{Preliminaries}\label{sec:preliminary}
Let $A \in \R^{m\times n}$ be a matrix of size $m\times n$ and has a rank of $\rho\leq \min(m, n)$. Let $k<\rho$ be a target rank to approximate $A$.    We  write the SVD decomposition of $A$ as
\[
A = U\left(
\begin{array}{cc}
\Sigma_1 &\mathbf 0 \\
\mathbf 0& \Sigma_2
\end{array}
\right)\left(
\begin{array}{c}
V_1^{\top} \\
V_2^{\top}
\end{array}
\right)
\]
where $\Sigma_1 \in \R^{k\times k}$, $\Sigma_2 \in \R^{(\rho-k)\times(\rho-k)}$, $V_1 \in \R^{n\times k}$ and $V_2 \in \R^{n\times (\rho-k)}$. We use $\sigma_1,\sigma_2,\ldots$ to denote the singular values of $A$ in the descending order, and $\lambda_{\max}(X)$ and $\lambda_{\min}(X)$ to denote the maximum and minimum eigen-values of a PSD matrix $X$. For any orthogonal matrix 
$U\in\R^{n\times \ell}$, let $U^{\perp}\in\R^{n\times (n-\ell)}$ denote an orthogonal matrix whose columns are an orthonormal basis spanning the subspace of $\R^n$ that is orthogonal to the column space of $U$.

 Let $\s = (s_1, \ldots, s_n)$ be a set of scores such that $\sum_{i=1}^ns_i=k$~\footnote{For the sake of discussion, we are not restricting the sum of these scores to be one but to be $k$, which does not affect our conclusions.}, one for each column of $A$. We will drawn $\ell$ independent samples with replacement from the set $[n]=\{1,\ldots, n\}$ using a multinomial distribution where the probability of choosing the $i$th column is $p_i=s_i/\sum_{j=1}^n s_j$. Let $i_1, \ldots, i_{\ell}$ be the indices of $\ell>k$ selected columns~\footnote{Note that some of the selected columns could be duplicate.}, and $S \in \R^{n\times \ell}$ be the corresponding sampling matrix, i.e, 
\begin{align*}
S_{i,j}=\left\{\begin{array}{cc}1, &\text{ if }i=i_j\\ 0,&\text{otherwise},
\end{array}\right.
\end{align*}
and $D\in\R^{\ell\times \ell}$ be a diagonal rescaling matrix with $\displaystyle D_{jj} =\frac{1}{\sqrt{s_{i_j}}}$.  
Given $S$, we construct the $C$ matrix as  
\begin{align}
C=AS=(A_{i_1}, \ldots, A_{i_\ell}).
\end{align}
Our interest is to bound the spectral norm error between $A$ and $P_CA$ for a given sampling matrix $S$, i.e., $
\|A - P_{C}A\|_2$, where $P_CA$ projects $A$ onto the column space of $C$.  For the benefit of presentation, we define  $\Omega = SD\in\R^{n\times\ell}$ to denote the sampling-and-rescaling matrix, and 
\begin{align}\label{eqn:Y}
Y = A\Omega,\;\quad \Omega_1 = V_1^{\top}\Omega, \;\quad \Omega_2 = V_2^{\top}\Omega, 
\end{align}
where $\Omega_1\in\R^{k\times \ell}$ and $\Omega_2\in\R^{(\rho-k)\times \ell}$.  Since the column space of $Y$ is the same to that of $C$, therefore 
\[
\|A - P_CA\|_2 = \|A - P_{Y}A\|_2
\]
and we will bound $ \|A - P_{Y}A\|_2$ in our analysis. 
%
Let $V^{\top}_1 = (\v_1,\ldots, \v_n)\in\R^{k\times n}$ and $V^{\top}_2=(\u_1,\ldots, \u_n)\in\R^{(\rho -k)\times n}$. 
It is easy to verify that
\begin{align*}
\Omega_1 = (\v_{i_1}, \ldots, \v_{i_{\ell}})D, \quad \Omega_2 = (\u_{i_1}, \ldots, \u_{i_{\ell}})D
\end{align*}
Finally, we let $\s^*=(s^*_1,\ldots, s^*_n)$ denote the SLS of $A$ relative to the best rank-$k$ approximation to $A$~\cite{mahoney-2011-randomized}, i.e.,  
$s_i^* = \|\v_i\|_2^2$. It is not difficult to show that $\sum_{i=1}^ns_i^*=k$. 

\section{Main Result}\label{sec:main}
Before presenting our main result,  we first characterize scores in $\s$ by two quantities as follows:
\begin{eqnarray}
c(\s) = \max\limits_{1 \leq i \leq n} \frac{s_i^*}{s_i}, \quad q(\s) = \max\limits_{1 \leq i \leq n} \frac{\sqrt{s_i^*}}{s_i} \label{eqn:cq}
\end{eqnarray}
Both quantities compare $\s$ to the SLS $\s^*$.  With  $c(\s)$ and $q(\s)$, we are ready to present our main theorem regarding the spectral error bound. 
\begin{thm} \label{thm:final-bound}
Let $A\in\R^{m\times n}$ has rank $\rho$ and  $C\in\R^{m\times \ell}$ contain the selected columns according to sampling scores in $\s$. With a probability $1 -\delta - 2k\exp(-\ell/[8kc(\s)])$, we have
\begin{align*}
&\|A - P_CA\|_2 \leq \sigma_{k+1}(1+\epsilon(\s))
\end{align*}
where $\epsilon(\s)$ is
\begin{align*}
\epsilon(\s)=3\left[\sqrt{c(\s)\frac{k(\rho+1-k)\log\left[\frac{\rho}{\delta}\right]}{\ell}} +q(\s)\frac{k\log\left[\frac{\rho}{\delta}\right]}{\ell}\right]\nonumber
\end{align*}
where $\sigma_{k+1} = \|A-A_k\|_2$ is the $(k+1)$th singular value of $A$. 
\end{thm}
\vspace*{-0.12in}
{\bf Remark:} Clearly, the spectral error bound and the successful probability in Theorem~\ref{thm:final-bound} depend on the quantities $c(\s)$ and $q(\s)$. In the subsection below, we study the two quantities to facilitate the understanding of the result in Theorem~\ref{thm:final-bound}. 

\subsection{More about the two quantities and their tradeoffs}
\vspace*{-0.05in}
The result in~Theorem~\ref{thm:final-bound} implies that the smaller the quantities $c(\s)$ and $q(\s)$, the better the error bound. Therefore, we first study  when  $c(\s)$ and $q(\s)$ achieve their minimum values. The key results are presented in the following two lemmas with their proofs deferred to the supplement. 
\begin{lemma}\label{lem:1}
The set of scores in $\s$ that minimize $q(\s)$ is given by $s_i\propto \sqrt{s^*_i}$, i.e., $s_i = \frac{k\sqrt{s^*_i}}{\sum_{i=1}^n\sqrt{s_i^*}}$.
\end{lemma}
{\bf Remark: } The sampling distribution with probabilities that are proportional to the square root of $s_i^*, i\in[n]$ falls in between the uniform sampling and the leverage-based sampling. 
\begin{lemma}\label{lem:2} $c(\s)\geq 1,\forall \s$ such that $\sum_{i=1}^ms_i=k$. The set of scores in $\s$ that minimize $c(\s)$ is given by $s_i =s_i^*$, and the minimum value of $c(\s)$ is $1$. 
\end{lemma}
Next, we discuss three special samplings with $\s$  (i) proportional to the square root of the SLS, i.e.,  $s_i\propto \sqrt{s_i^*}$ (referred to as square-root leverage-based sampling or {\bf sqL-sampling} for short), (ii) equal to the SLS, i.e., $s_i = s_i^*$ (referred to as leverage-based sampling or {\bf L-sampling} for short), and (iii) equal to uniform scalars $s_i = k/n$ (referred to as uniform sampling or \textbf{U-sampling} for short).  Firstly, if $s_i\propto \sqrt{s_i^*}$ , $q(\s)$ achieves its minimum value and we have the two quantities written as
\begin{equation}\label{eqn:sq}
\begin{split}
q_{sqL}&=  \frac{1}{k}\sum_{i=1}^n\sqrt{s_i^*}\\
c_{sqL} &= \max_{i}\frac{s_i^*\sum_i \sqrt{s_i^*}}{k\sqrt{s_i^*}}=q_{sqL}\max_i\sqrt{s_i^*}\\
\end{split}
\end{equation}
In this case, when $\s^*$ is flat (all SLS are equal), then $q_{sqL} =\sqrt{\frac{n}{k}}$ and $c^{sqL}=1$. The bound becomes $\widetilde O(\sqrt{(\rho+1-k)k/\ell} + \sqrt{nk/\ell^2})\sigma_{k+1}$ that suppresses  logarithmic terms. To analyze $q_{sqL}$ and $c_{sqL}$ for skewed SLS, we consider a power-law distributed SLS, i.e., there exists a small constant $a$ and power index $p >2$, such that $s^*_{[i]}, i=1,\ldots, n$ ranked in descending order satisfy 
\[
s^*_{[i]}\leq a^2i^{-p}, \quad i=1,\ldots, n
\]
Then it is not difficult to show that 
\[
\frac{1}{k}\sum_{i=1}^n\sqrt{s_i}\leq \frac{a}{k}\left(1 + \frac{2}{p-2}\right)
\]
which is independent of $n$. Then the error bound in Theorem~\ref{thm:final-bound} becomes $O\left(\sqrt{\frac{\rho+1-k}{\ell}} + \frac{1}{\ell}\right)\sigma_{k+1}$, which is better than that in~(\ref{eqn:bound-2}).

Secondly, if $s_i\propto s_i^*$, then  $c(\s)$ achieves its minimum value and we have the two quantities written as 
\begin{equation}\label{eqn:sls}
q_{L} = \max_i \frac{1}{\sqrt{s_i^*}}, \quad c_L= 1
\end{equation}
In this case, when $\s^*$ is flat, we have $q_L= \sqrt{\frac{n}{k}}$ and $c_L =1$ and  the same bound  $\widetilde O(\sqrt{(\rho+1-k)k/\ell} + \sqrt{nk/\ell^2})\sigma_{k+1}$  follows. However, when $\s^*$ is skewed, i.e., there exist very small SLS, then $q_L$ could be very large.  As a comparison, the $q(\s)$ for  sqL-sampling is always smaller than that for L-sampling due the following inequality 
\begin{align*}
q_{sqL}&=\frac{1}{k}\sum_{i=1}^n\sqrt{s_i^*} = \frac{1}{k}\sum_{i=1}^n\frac{s_i^*}{\sqrt{s_i^*}}< \max_i \frac{1}{\sqrt{s_i^*}}\frac{\sum_{i=1}^ns_i^*}{k}\\
&=\max_i \frac{1}{\sqrt{s_i^*}} = q_{L}
\end{align*}

\begin{figure}[t]
\centering
\includegraphics[scale=0.5]{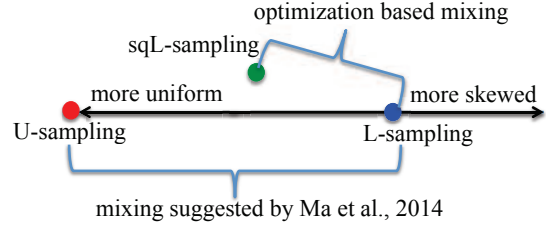}
\caption{An illustration of different sampling strategies. The mixing strategy suggested by~\cite{DBLP:conf/icml/MaMY14} is a convex combination of U-sampling and L-sampling. Our optimization approach gives an intermediate sampling between the sqL-sampling and the L-sampling. }\label{fig:sample}
\vspace*{-0.3in}
\end{figure}

Lastly,  we consider the uniform sampling $s_i=\frac{k}{n}$ . Then the two quantities become
\begin{equation}\label{eqn:unif}
q_{U}= \max_i\frac{n\sqrt{s_i^*}}{k},\quad c_{U} = \max_i\frac{ns_i^*}{k}
\end{equation}
Similarly, if $\s_*$ is flat,  $q_U =\sqrt{\frac{n}{k}}$ and $c_U = 1$.  Moreover, it is interesting to compare the two quantities for the sqL-sampling in~(\ref{eqn:sq}) and for the uniform sampling in~(\ref{eqn:unif}). 
\begin{align*}
q_{sqL}&=  \frac{1}{k}\sum_{i=1}^n\sqrt{s_i^*}\leq \max_i\frac{n\sqrt{s_i}}{k} = q_U\\
c_{sqL}& =  \max_i \frac{1}{k}\sqrt{s_i^*}\sum_{i=1}^n\sqrt{s_i^*}\leq \max_i \frac{ns_i^*}{k} = c_U
\end{align*}
From the above discussions, we can see that when $\s_*$ is a flat vector, there is no difference between the three sampling scores for $\s$. The difference comes from when $\s_*$ tends to be skewed. In this case, $s_i\propto \sqrt{s_i^*}$ works almost for sure better than uniform distribution and could also be potentially better than $s_i \propto s_i^*$  according to the sampling dependent error bound in Theorem~\ref{thm:final-bound}. A similar tradeoff between the L-sampling and U-sampling but with  a different taste was observed  in~\cite{DBLP:conf/icml/MaMY14}, where they showed that  for least square approximation by CSS  leveraging-based least square estimator could have a  large variance when there exist very small SLS. Nonetheless, our bound here exhibits more insights, especially on the sqL-sampling.  More importantly, the sampling dependent bound renders the flexibility in choosing the sampling scores by adjusting them according to the distribution of the SLS. In next subsection, we present an optimization approach to find better sampling scores. In Figure~\ref{fig:sample}, we give a quick view of different sampling strategies.

\subsection{Optimizing the error bound}
\vspace*{-0.05in}
As indicated by the result in Theorem~\ref{thm:final-bound}, in order to achieve a good performance, we need to make a balance between $c(\s)$ an $q(\s)$, where $c(\s)$ affects  not only the error bound but also the successful probability. To address this issue, we propose a constrained optimization approach.  More specifically, to ensure that the failure probability is no more than $3\delta$, we impose the following constraint on $c(\s)$
\begin{align}\label{eqn:constraint}
&\frac{\ell}{8kc(\s)}\geq \log\left(\frac{k}{\delta}\right),\: i.e.,\:\max_i \frac{s_i^*}{s_i}\leq \frac{\ell}{8k \log\left(\frac{k}{\delta}\right) } :=\gamma
\end{align}
Then we cast the problem into minimizing $q(\s)$ under the constraint in (\ref{eqn:constraint}), i.e.,
\begin{align}
&\min\limits_{\s \in \R_+^n} \; \max\limits_{1 \leq i \leq n} \frac{\sqrt{s_i^*}}{s_i}\nonumber\\
&\mbox{s.t.} \quad \s^{\top}\mathbf{1} = k, \; s_i^* \leq \gamma s_i, i=1, \ldots, n \label{eqn:opt}
\end{align}
It is easy to verify that the optimization problem in (\ref{eqn:opt}) is convex. Next, we develop an efficient bisection search algorithm to solve the above problem with a linear convergence rate. To this end, we introduce a slack variable $t$ and rewrite the optimization problem in (\ref{eqn:opt}) as
\begin{equation} \label{eqn:opt-1}
\begin{aligned} 
\min\limits_{\s \in \R_+^n, t \geq 0}&\quad  t, \quad \mbox{s.t.}\quad \s^{\top}\mathbf{1} = k \\
\mbox{and}&\quad  \frac{s_i^*}{s_i} \leq \min\left(\gamma, t\sqrt{s_i^*}\right), i=1, \ldots, n
\end{aligned}
\end{equation}
We now find the optimal solution by performing bisection search on $t$. Let $t_{\max}$ and $t_{\min}$ be the upper and lower bounds for $t$. We set $t = (t_{\min} + t_{\max})/2$ and decide the feasibility of $t$ by simply computing the quantity
\[
 f(t) = \sum_{i=1}^n \frac{s_i^*}{\min\left(\gamma, t\sqrt{s_i^*}\right)}
\]
Evidently, $t$ is a feasible solution if $f(t) \leq k$ and is not if $f(t) > k$. Hence, we will update $t_{\max} = t$ if $f(t) \leq k$ and $t_{\min} = t$ if $f(t) > k$. To run the bisection algorithm, we need to decide  initial $t_{\min}$ and $t_{\max}$. We can set $t_{\min} = 0$. To compute $t_{\max}$, we make an explicit construction of $\s$ by distributing the $(1 - \gamma^{-1})$ share of the largest element of $\s_*$ to the rest of the list. More specifically, let $j$ be the index for the largest entry in $\s^*$. We set $s_{j} = \|\s^*\|_{\infty}\gamma^{-1}$ and $s_i = s_i^* + (1 - \gamma^{-1})\|\s^*\|_{\infty}/(n - 1)$ for $i \neq j$. Evidently, this solution satisfies the constraints $s_i^* \leq \gamma s_i, i\in[n]$ for $\gamma \geq 1$. With this construction, we can show that 
\begin{align*}
q(\s)  \leq \max\left(\frac{\gamma}{\sqrt{\|\s^*\|_\infty}}, \frac{n - 1}{\sqrt{\|\s^*\|_{\infty}}(1 - \gamma^{-1})}\right)
\end{align*}
Therefore, we set initial $t_{\max}$ to the value in R.H.S of the above inequality.  Given the optimal value of $t=t_*$ we compute  the optimal value of $s_i$ by $s_i=\frac{s_i^*}{\min(\gamma, t_*\sqrt{s_i^*})}.$ The corresponding sampling distribution clearly lies between L-sampling  and sqL-sampling. In particular, when $\gamma=1$ the resulting sampling distribution is L-sampling due to Lemma~\ref{lem:2} and when $\gamma\rightarrow\infty$ the resulting sampling distribution approaches sqL-sampling. 

Finally, we comment on the value of $\ell$. In order to make the constraint in~(\ref{eqn:constraint}) feasible, we need to ensure $\gamma\geq 1$. Therefore, we need $\ell\geq \Omega(k\log\left(\frac{k}{\delta}\right))$.

\subsection{Subsequent Applications}
\vspace*{-0.05in}
Next, we discuss  two subsequent applications of CSS, one for low rank approximation and one for least square approximation.  

\textbf{Rank-$k$ approximation.}  If a rank-$k$ approximation is desired, we need to do some postprocessing since $P_CA$ might has rank larger than $k$. We can use the same algorithm as presented in \cite{conf/focs/BoutsidisDM11}. In particular, given the constructed $C\in\R^{n\times\ell}$, we first orthonormalize the columns of $C$ to construct a matrix $Q\in\R^{m\times \ell}$ with orthonormal columns, then compute the best rank-$k$ approximation of $Q^{\top}A\in\R^{\ell\times n}$ denoted by $(Q^{\top}A)_k$, and finally construct the low-rank approximation as $Q(Q^{\top}A)_k$. It was shown that (Lemma 2.3 in~\citep{conf/focs/BoutsidisDM11})
\[
\|A - Q(Q^{\top}A)_k\|_2\leq \sqrt{2}\|A - \Pi^2_{C,k}(A)\|_2
\]
where $\Pi^2_{C,k}(A)$ is the best approximation to $A$ within the column space of $C$ that has rank at most $k$. The running time of above procedure is $O(mn\ell + (m+n)\ell^2)$. Regarding its error bound,  the above inequality together with the following theorem  implies that its spectral error bound  is only amplified by a factor of $\sqrt{2}$ compared to that of $P_CA$. 
\begin{thm} \label{thm:final-bound-2}
Let $A\in\R^{m\times n}$ has rank $\rho$  and  $C\in\R^{m\times \ell}$ contain the selected columns according to sampling scores in $\s$. With a probability $1 - \delta - 2k\exp(-\ell/[8kc(\s)])$, we have
\begin{align*}
&\|A - \Pi^2_{C,k}(A)\|_2 \leq\sigma_{k+1}(1+\epsilon(\s))\nonumber
\end{align*}
where $\epsilon(\s)$ is given in Theorem~\ref{thm:final-bound}.  
\end{thm}
\textbf{Least Square Approximation.}
CSS has been used in least square approximation for developing  faster and interpretable algorithms. In these applications, an over-constrained least square problem is considered, i.e., given $A\in\R^{m\times n}$ and $\b\in\R^m$ with $m\gg n$, to solve the following problem:
\begin{align}\label{eqn:ls}
\x_{opt}=\arg\min_{\x\in\R^n}\|A\x - \b\|^2_2
\end{align}
The procedure for applying CSS to least square approximation is (i) to sample a set of $\ell>n$ rows from $A$ and form a sampling-and-rescaling matrix denoted by $\Omega\in\R^{\ell\times m}$~\footnote{We abuse the same notation $\Omega$.}; (ii) to solve the following reduced least square problem: 
\begin{align}\label{eqn:lsapprox}
\xh_{opt}=\arg\min_{\x\in\R^n}\|\Omega A\x - \Omega\b\|^2_2
\end{align}
It is worth pointing out  that in this case the SLS $\s^*=(s^*_1,\ldots, s_m^*)$ are computed based on the the left singular vectors $U$ of $A$ by  $s_i^* =\|U_{i*}\|_2^2$, where $U_{i*}$ is the $i$-th row of $U$. One might be  interested to see whether we can apply our analysis to derive a sampling dependent error bound for the approximation error $\|\x_{opt} - \xh_{opt}\|_2$ similar to previous bounds of the form $\|\x_{opt} - \xh_{opt}\|_2\leq \frac{\epsilon}{\sigma_{min}(A)}\|A\x_{top}-\b\|_2$. Unfortunately, naively combining our analysis with previous analysis is a worse case analysis, and consequentially yields a worse bound. The reason will become clear in our later discussions. However, the statistical analysis in~\cite{DBLP:conf/icml/MaMY14} does indicate that $\xh_{opt}$ by using sqL-sampling could have smaller variance than that using L-sampling.

\section{Numerical Experiments}\label{sec:empirical}
Before delving into the detailed analysis, we present some experimental results. We consider synthetic data with the data matrix $A$ generated from one of the three different classes of distributions introduced below, allowing the SLS vary from nearly uniform to very nonuniform.
\vspace*{-0.1in} 
\begin{itemize}
\item Nearly uniform SLS (GA). Columns of $A$ are generated from a multivariate normal distribution $\mathcal N(\mathbf 1_m, \Sigma)$, where $\Sigma_{ij}=2*0.5^{|i-j|}$. This data is referred to as GA data. 
\item Moderately nonuniform SLS ($T_3$). Columns of $A$ are generated from a multivariate $t$-distribution with $3$ degree of freedom and covariance matrix $\Sigma$ as before. This data is referred to as $T_3$ data. 
\item Very nonuniform SLS ($T_1$). Columns of $A$ are generated from a multivariate $t$-distribution with $1$ degree of freedom and covariance matrix $\Sigma$ as before. This data is referred to as $T_1$ data. 
\end{itemize}
\vspace*{-0.15in}
These distributions  have been used in~\cite{DBLP:conf/icml/MaMY14} to generate synthetic data for empirical evaluations. 
\begin{figure}[t]
\includegraphics[scale=0.17]{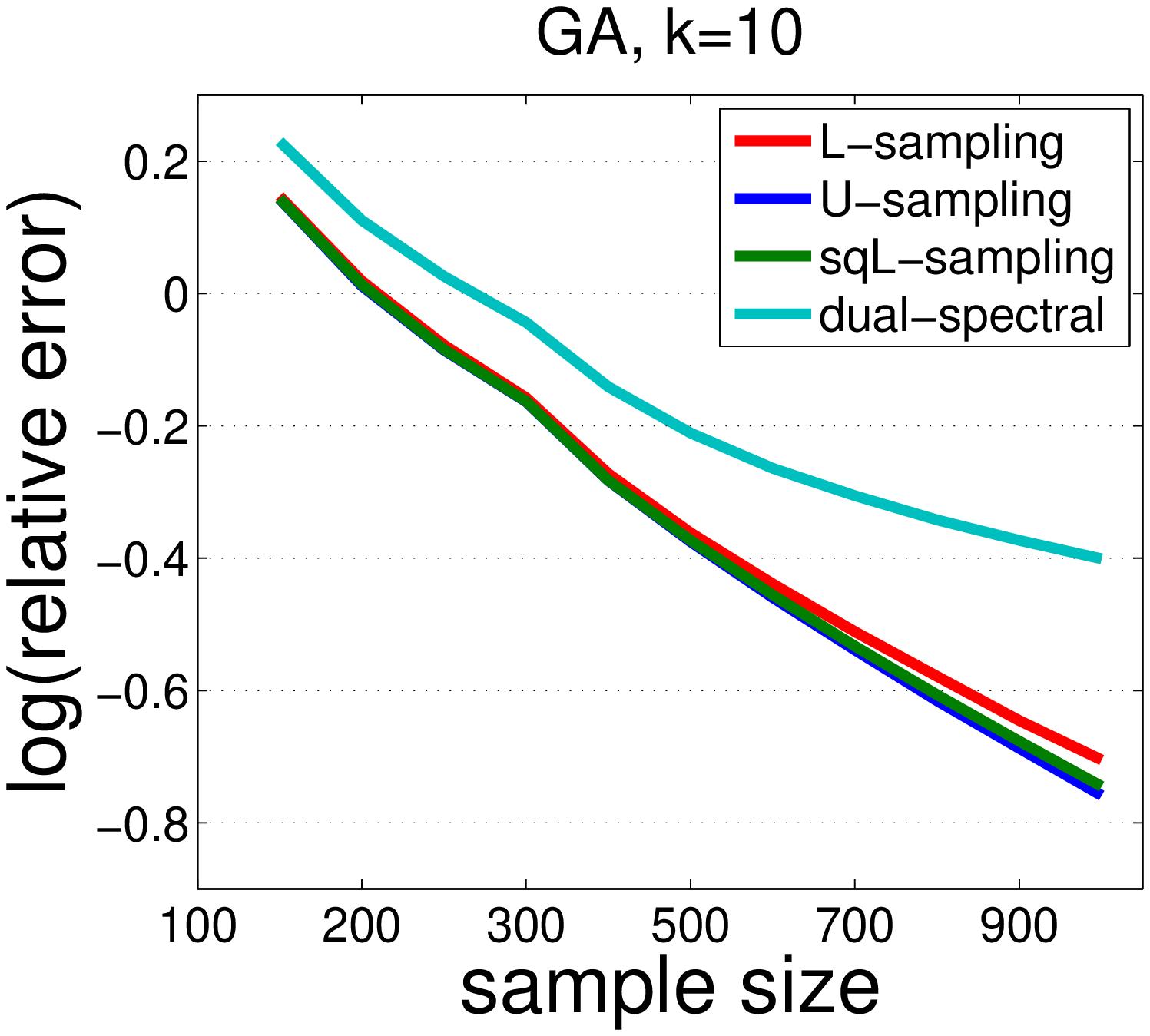}\hspace*{-0.1in}
\includegraphics[scale=0.17]{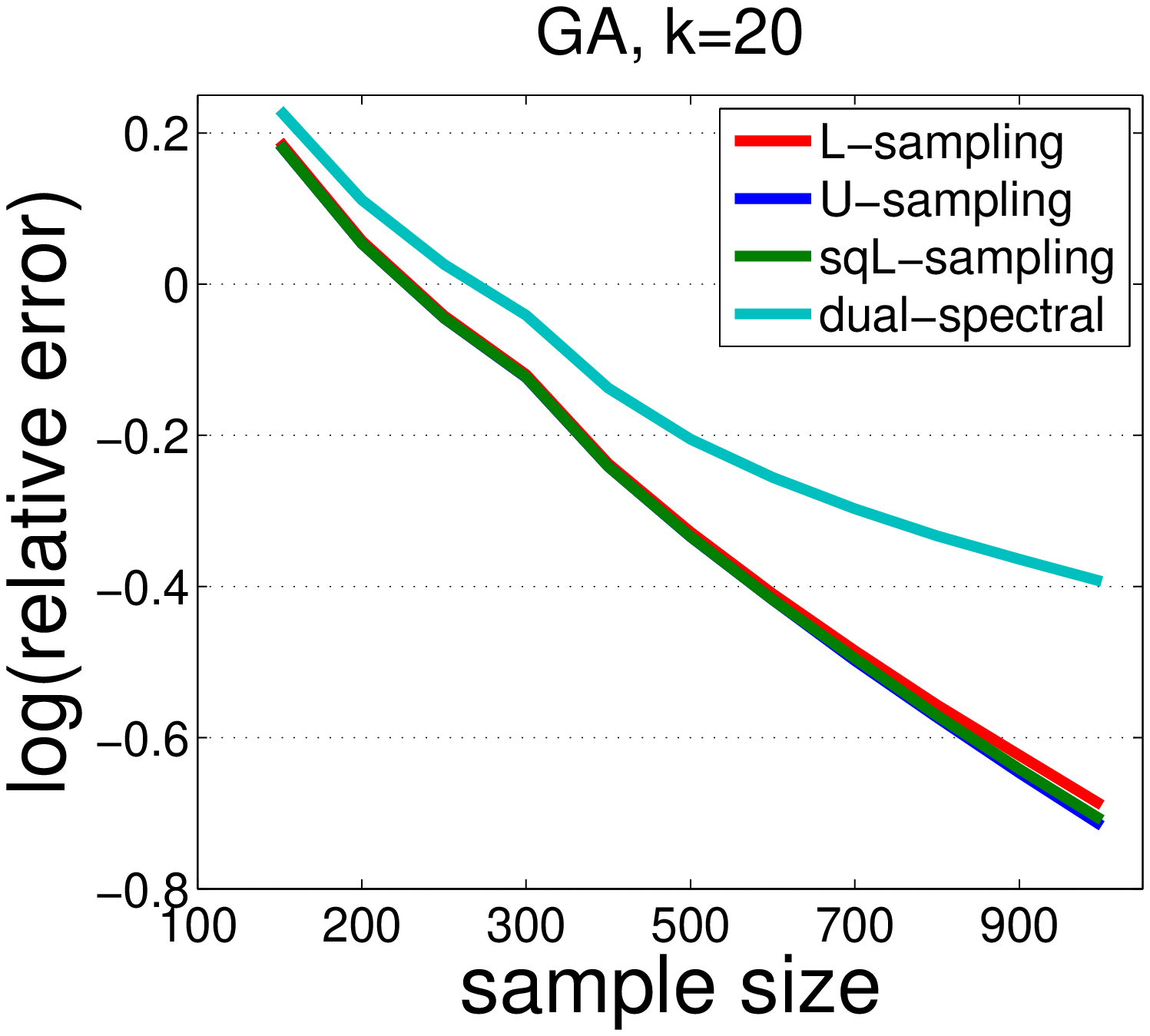}\hspace*{-0.1in}
\includegraphics[scale=0.17]{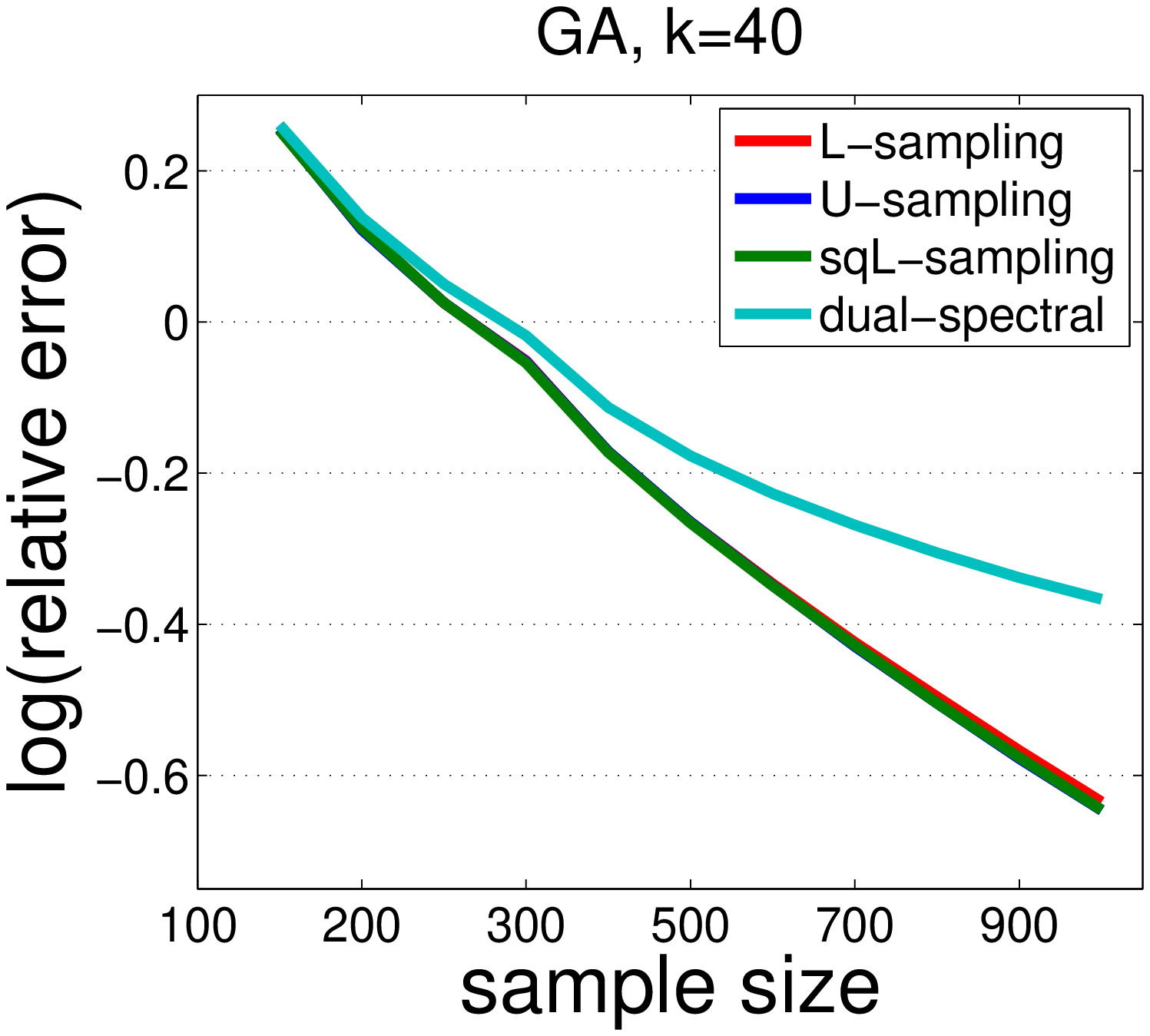}

\includegraphics[scale=0.17]{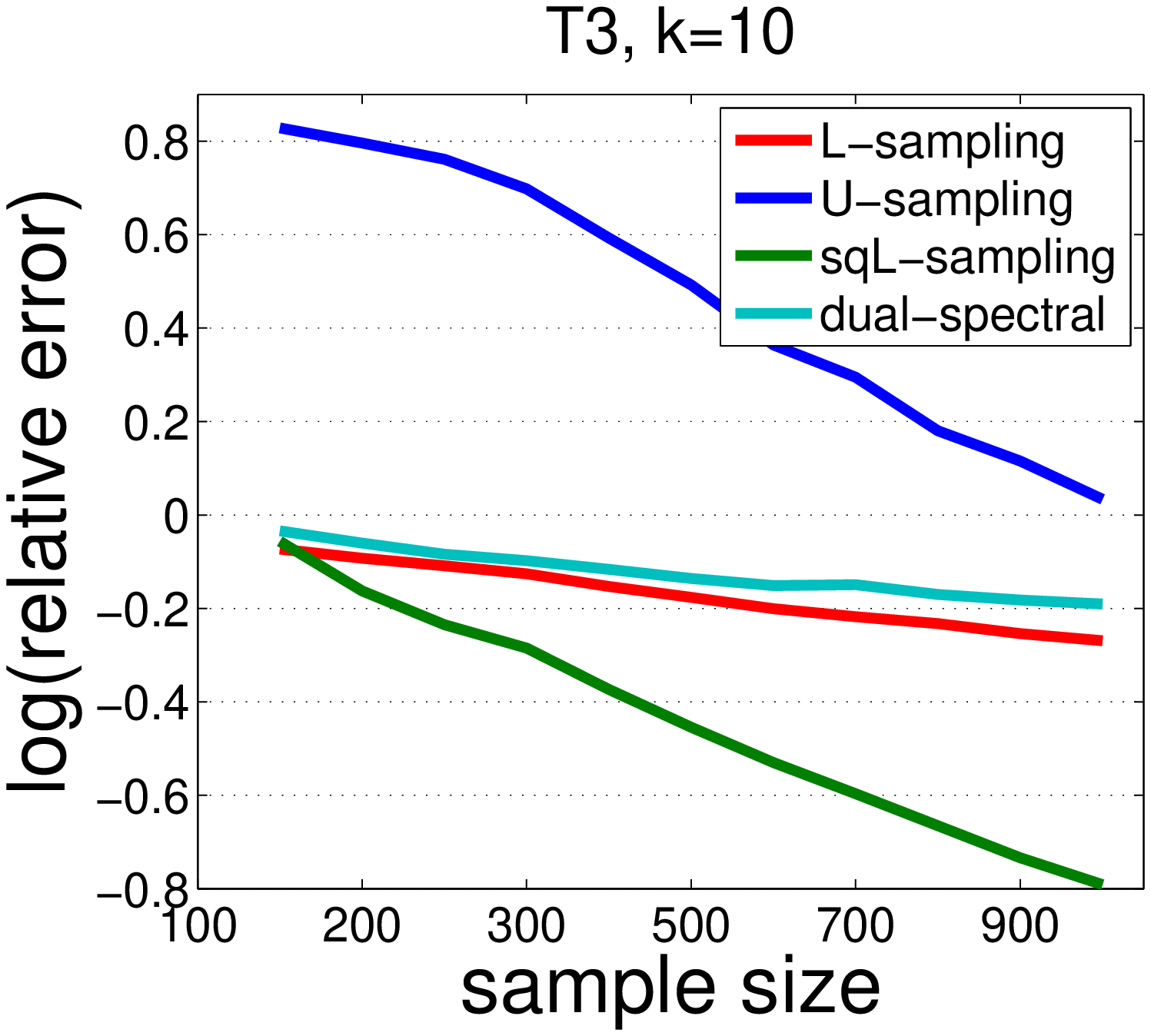}\hspace*{-0.1in}
\includegraphics[scale=0.17]{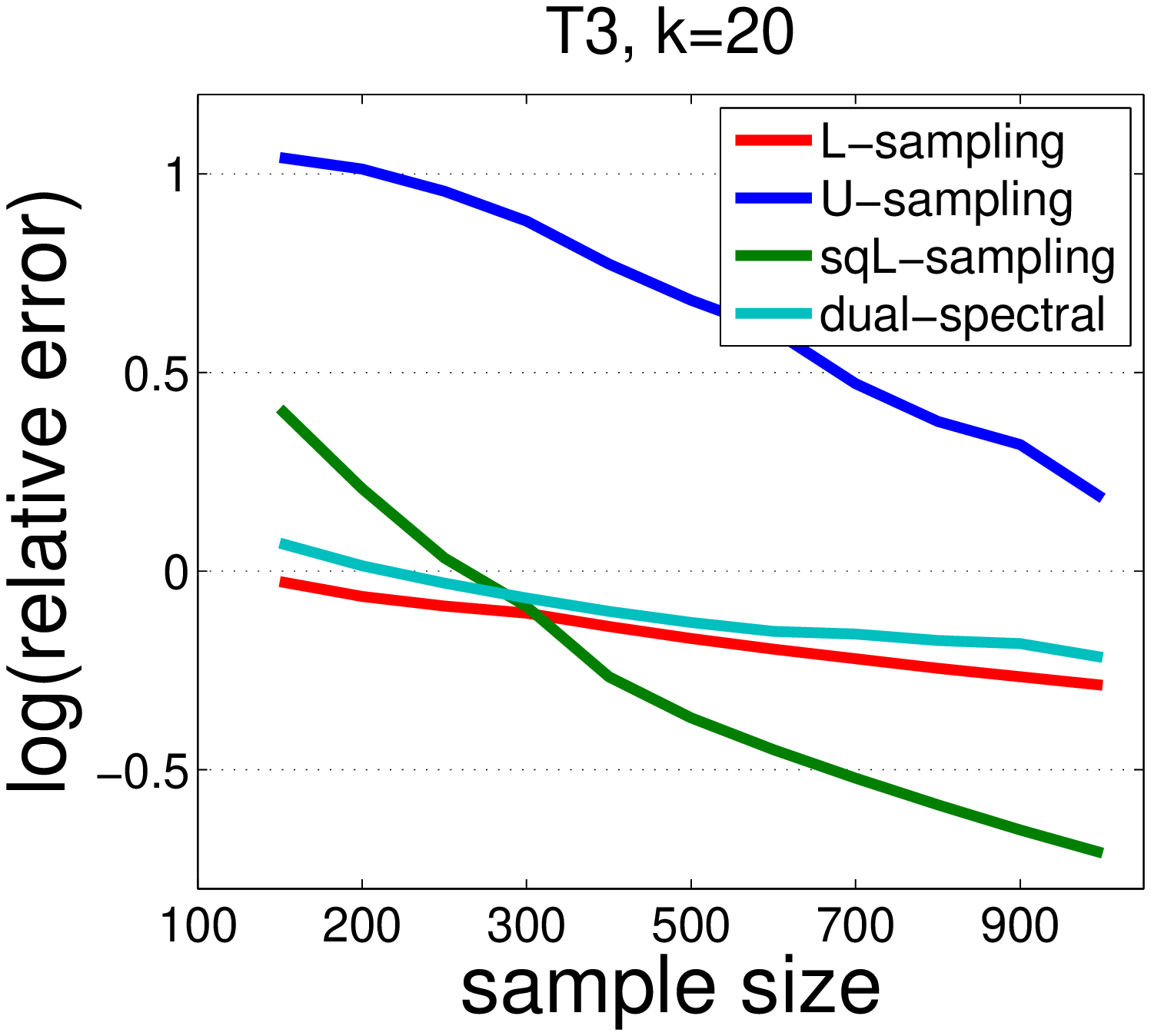}\hspace*{-0.1in}
\includegraphics[scale=0.17]{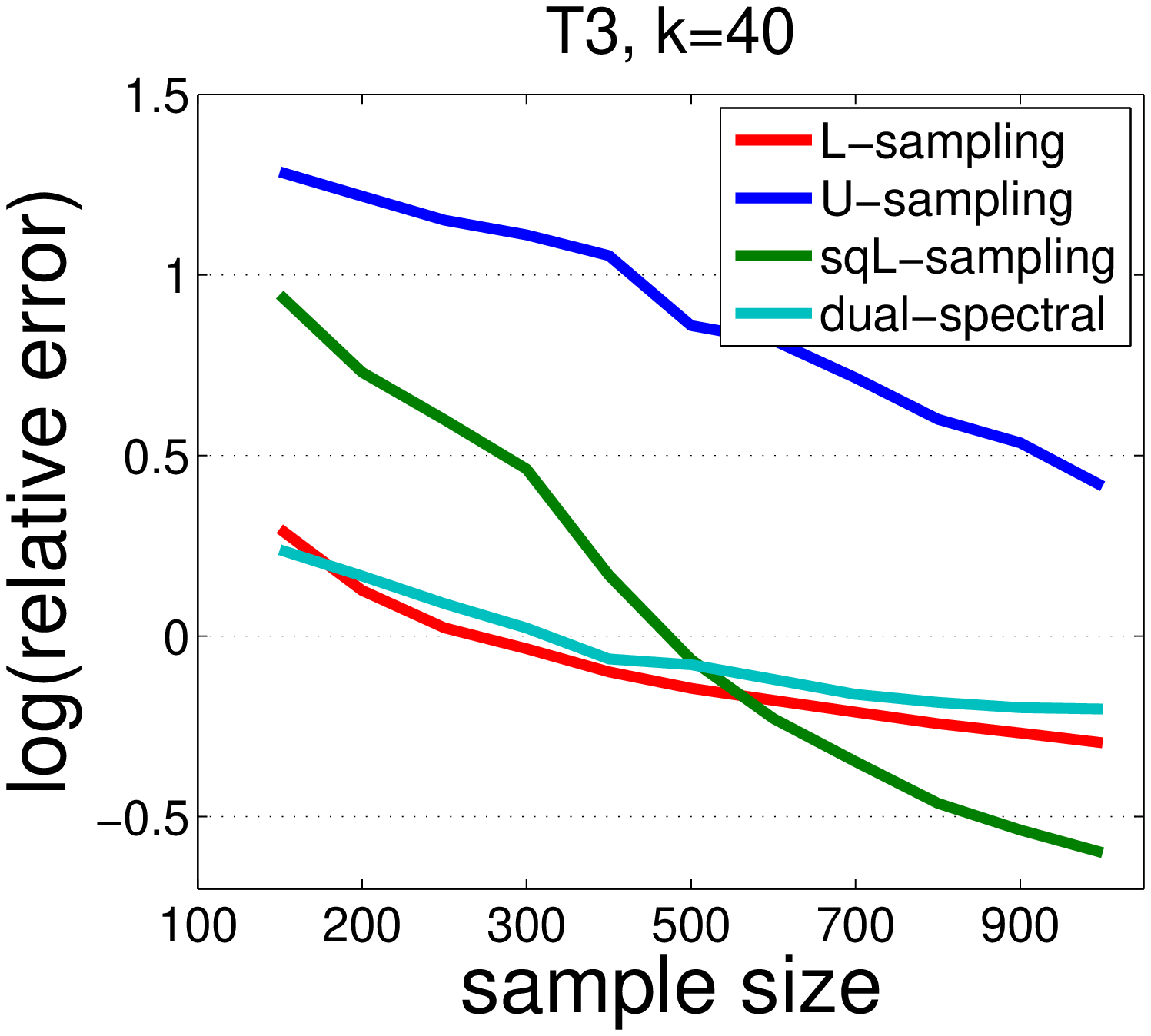}

\includegraphics[scale=0.17]{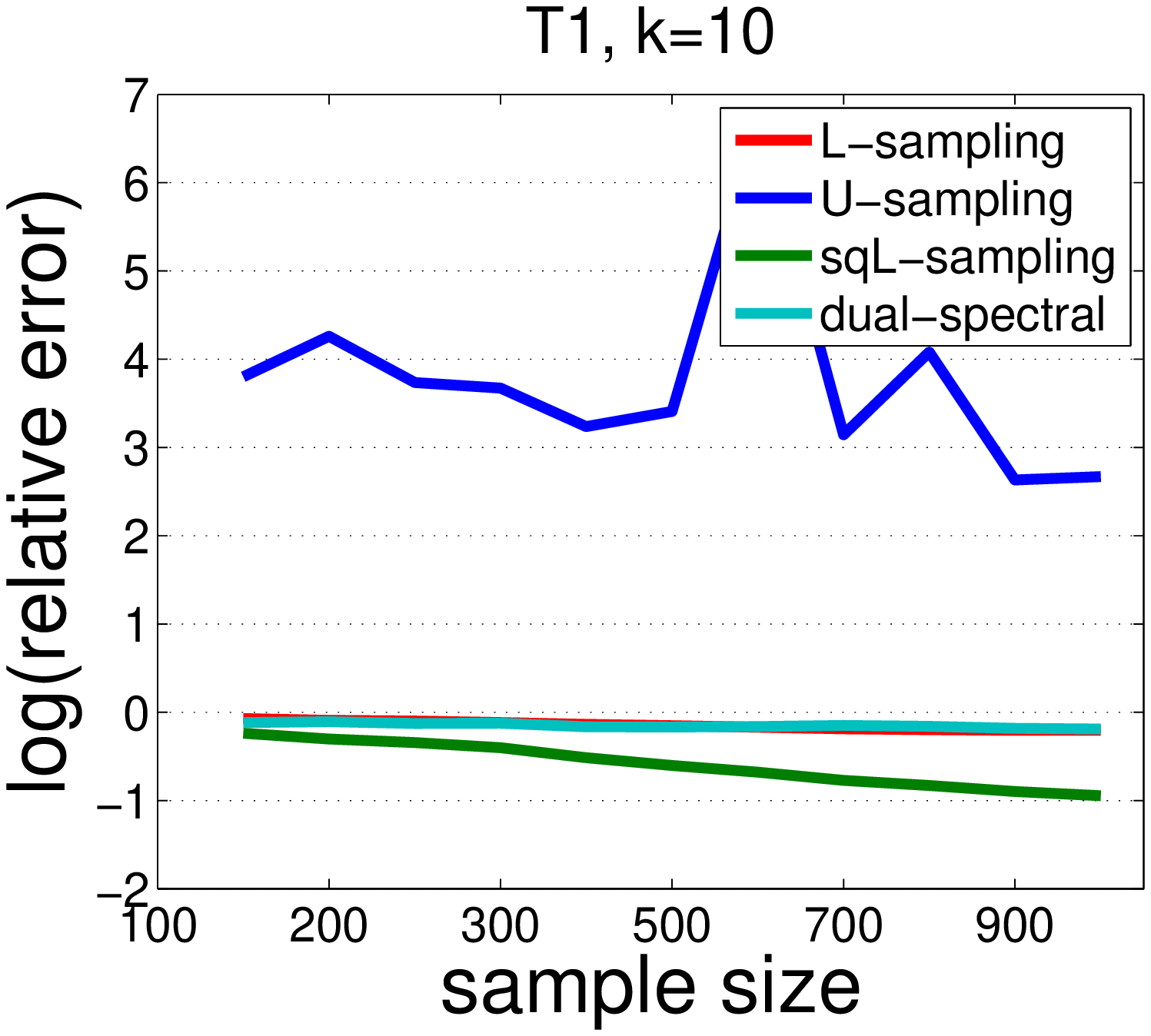}\hspace*{-0.1in}
\includegraphics[scale=0.17]{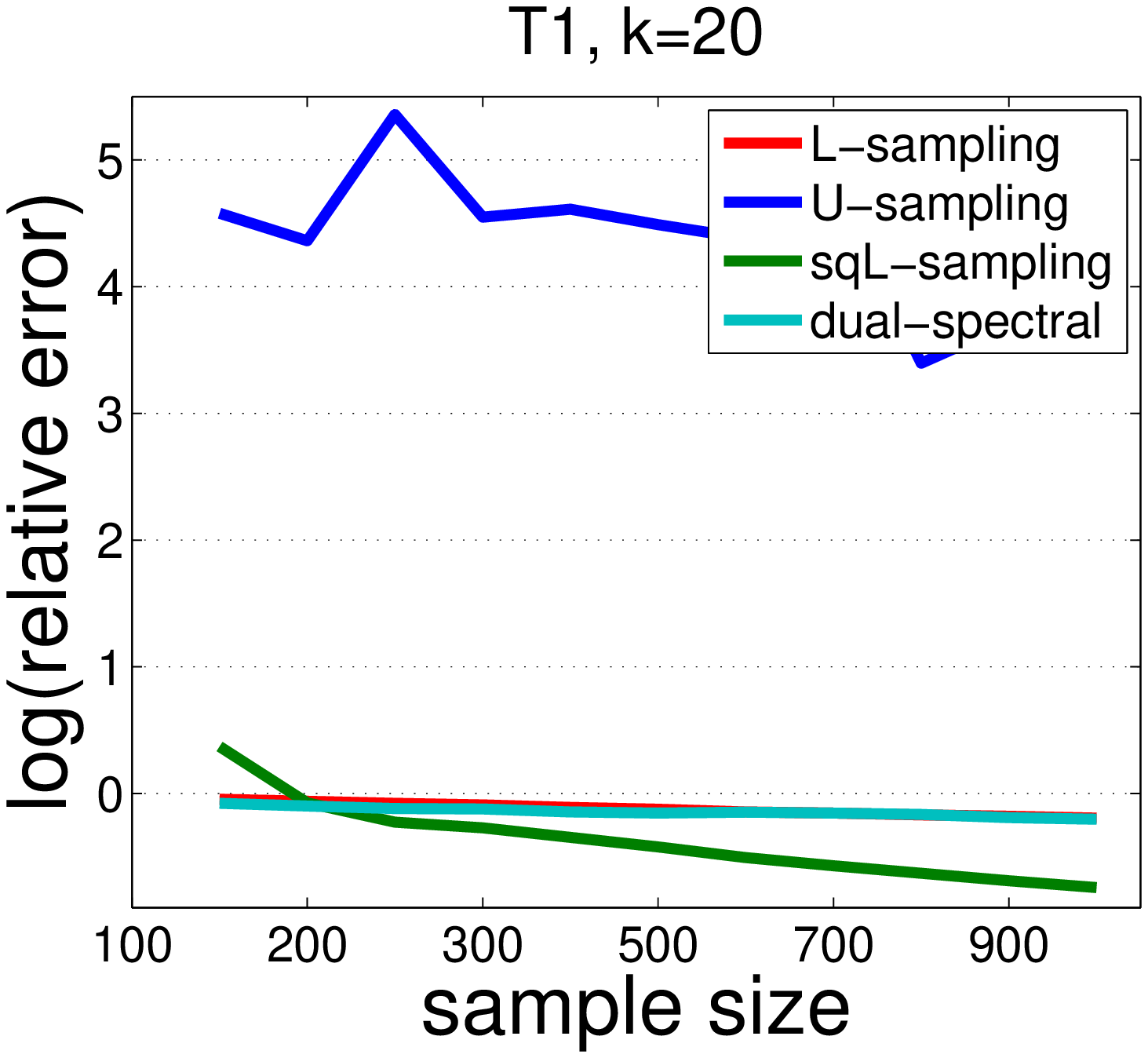}\hspace*{-0.1in}
\includegraphics[scale=0.17]{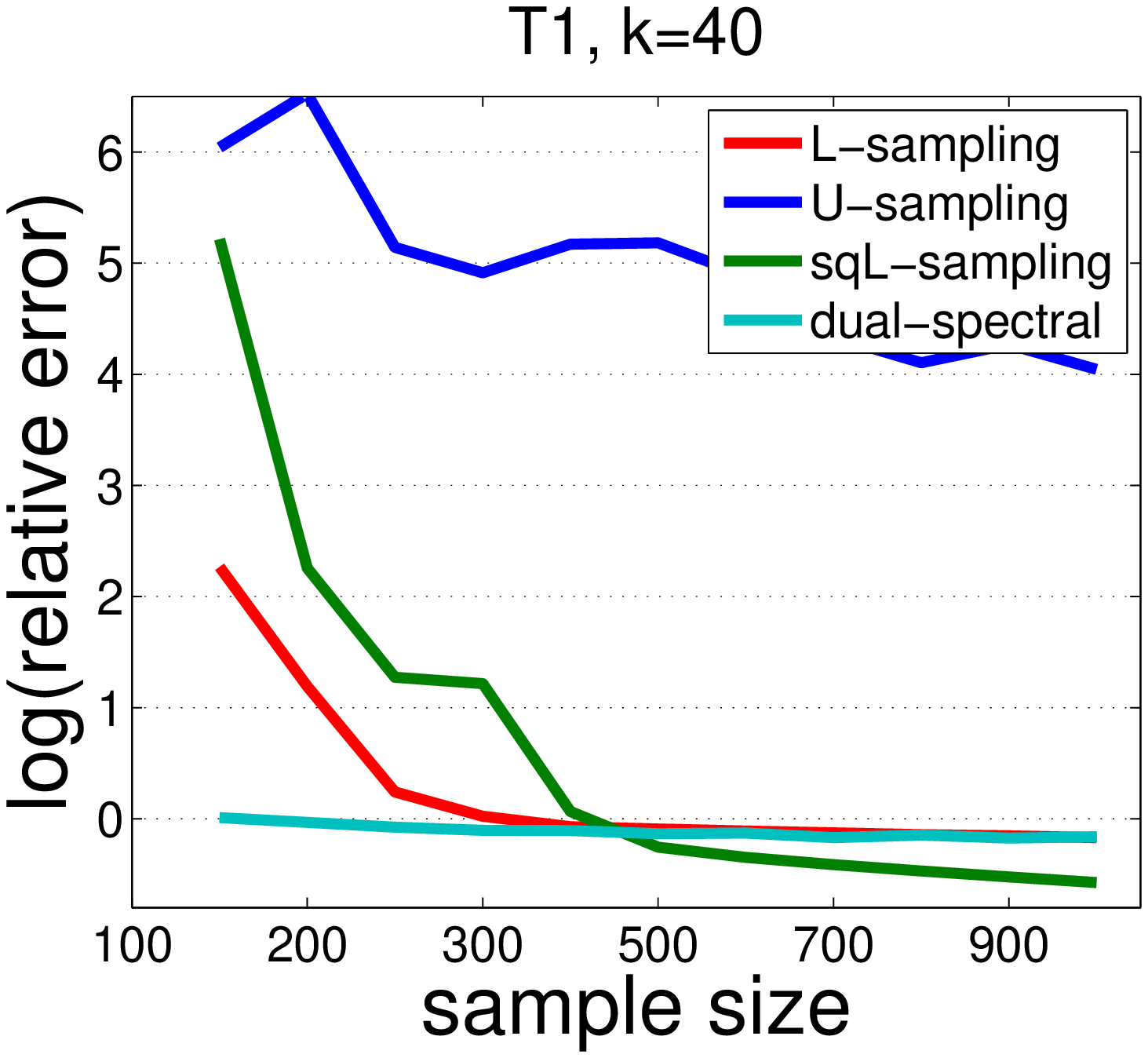}

\caption{Comparison of the spectral error for different data, different samplings, different target rank and different sample size.}\label{fig:1}
\vspace*{-0.1in}
\end{figure}

\begin{figure}[t]
\centering
\includegraphics[scale=0.17]{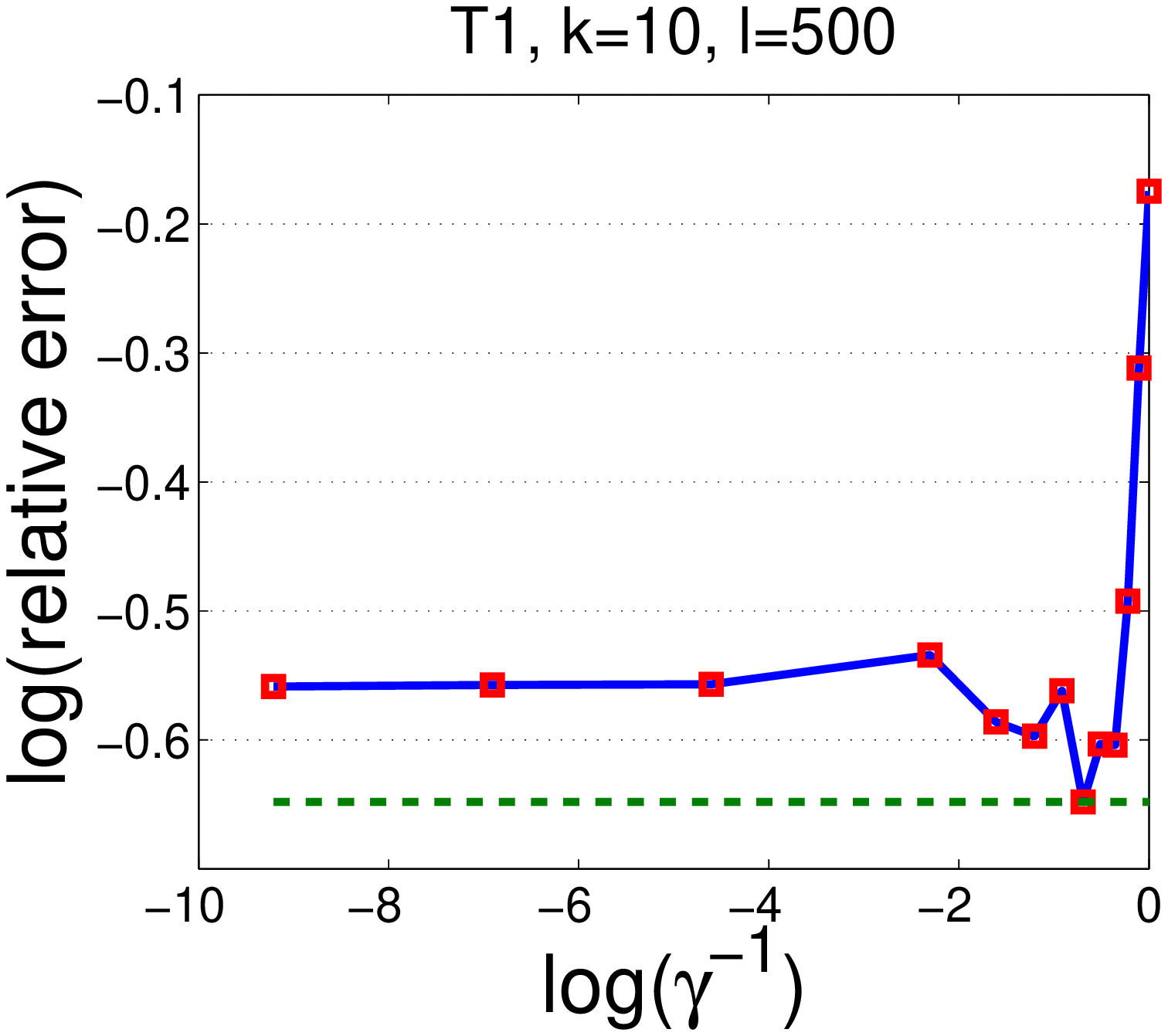}\hspace*{-0.1in}
\includegraphics[scale=0.17]{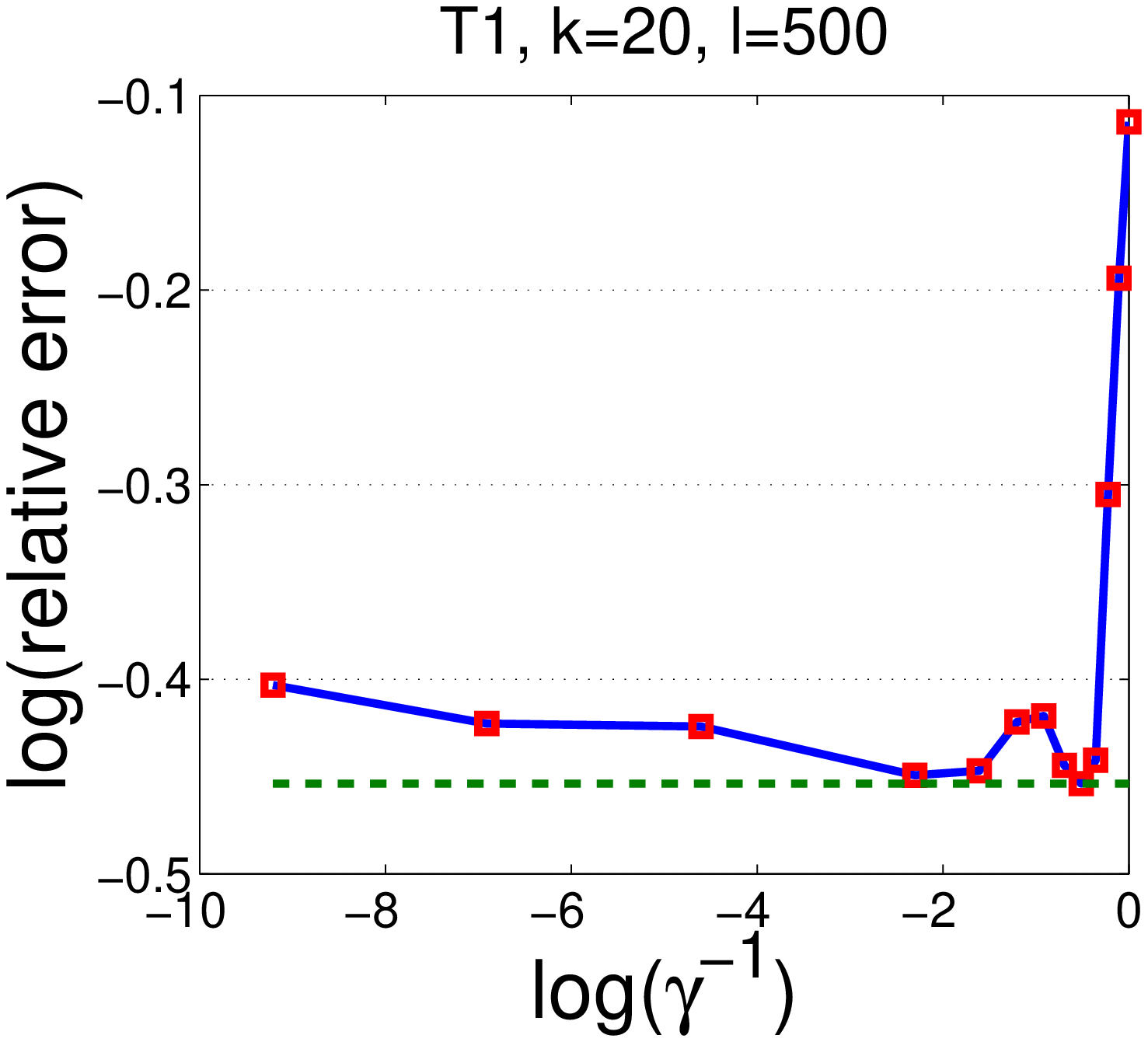}\hspace*{-0.1in}
\includegraphics[scale=0.17]{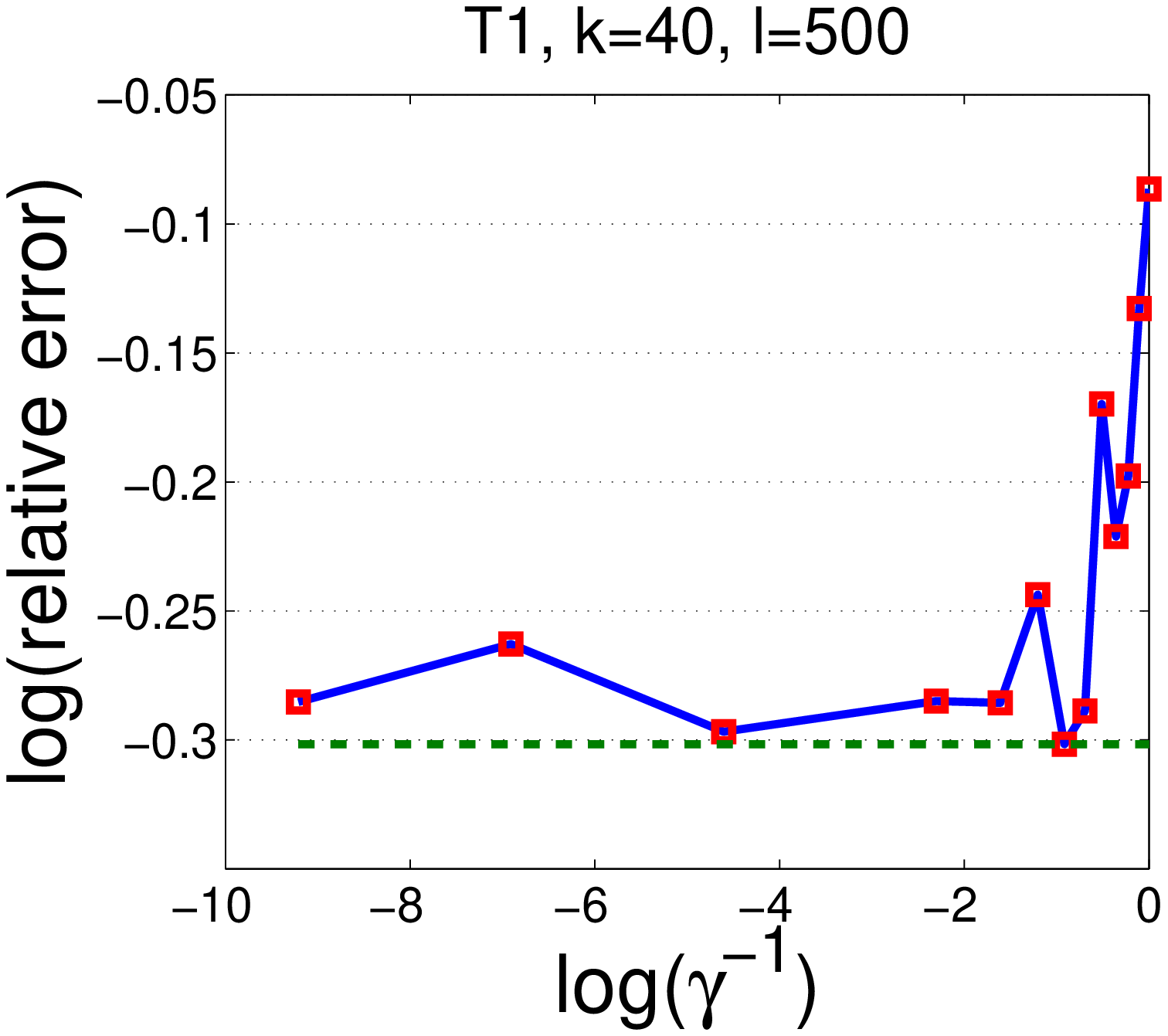}

\caption{The spectral error for the sampling probabilities found by the constrained optimization approach with different values of $\gamma\geq 1$. The left most point corresponds to sqL-sampling and the right most point corresponds to L-sampling.}\label{fig:2}
\vspace*{-0.2in}
\end{figure}
We first compare the spectral norm reconstruction error of the three different samplings, namely L-sampling, U-sampling and the sqL-sampling, and the deterministic dual set spectral sparsification algorithm. We generate synthetic data with $n=m=1000$ and repeat the experiments 1000 times. We note that the rank of the generated data matrix is $1000$. The averaged results are shown in  Figure~\ref{fig:1}.  From these results we observe  that (i) when the SLS are nearly uniform, the three sampling strategies perform similarly as expected; (ii) when the SLS become nonuniform, sqL-sampling performs always better than U-sampling and better than the L-sampling when the target rank is small (e.g., $k=10$) or the sample size $\ell$ is large; (iii) when the SLS are non-uniform, the spectral norm reconstruction error of sqL-sampling decreases faster than L-sampling w.r.t the sample size $\ell$; (iv) randomized algorithms generally perform  better than the deterministic dual set sparsification algorithm. 
\begin{figure}[t]
\includegraphics[scale=0.17]{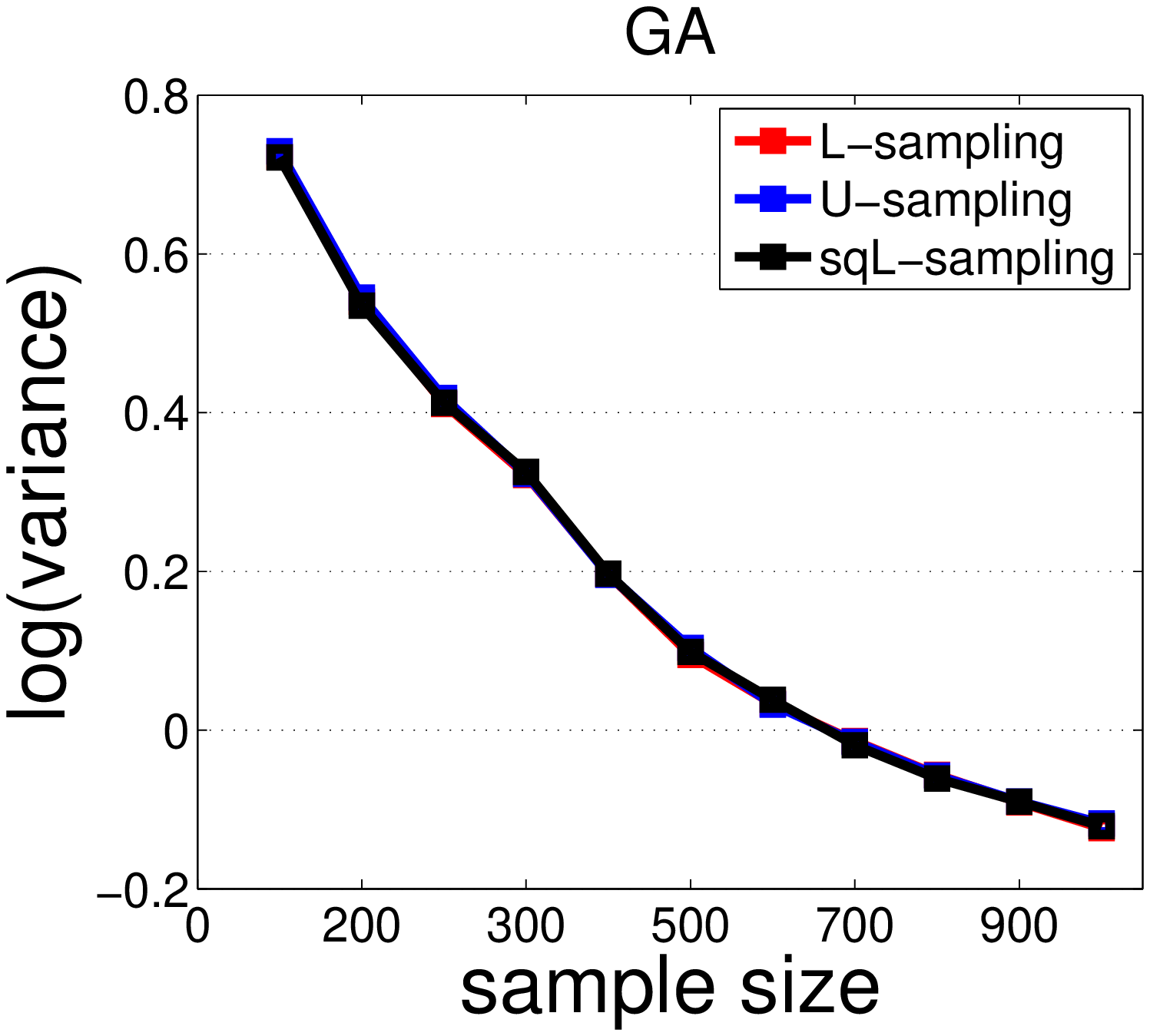}\hspace*{-0.1in}
\includegraphics[scale=0.17]{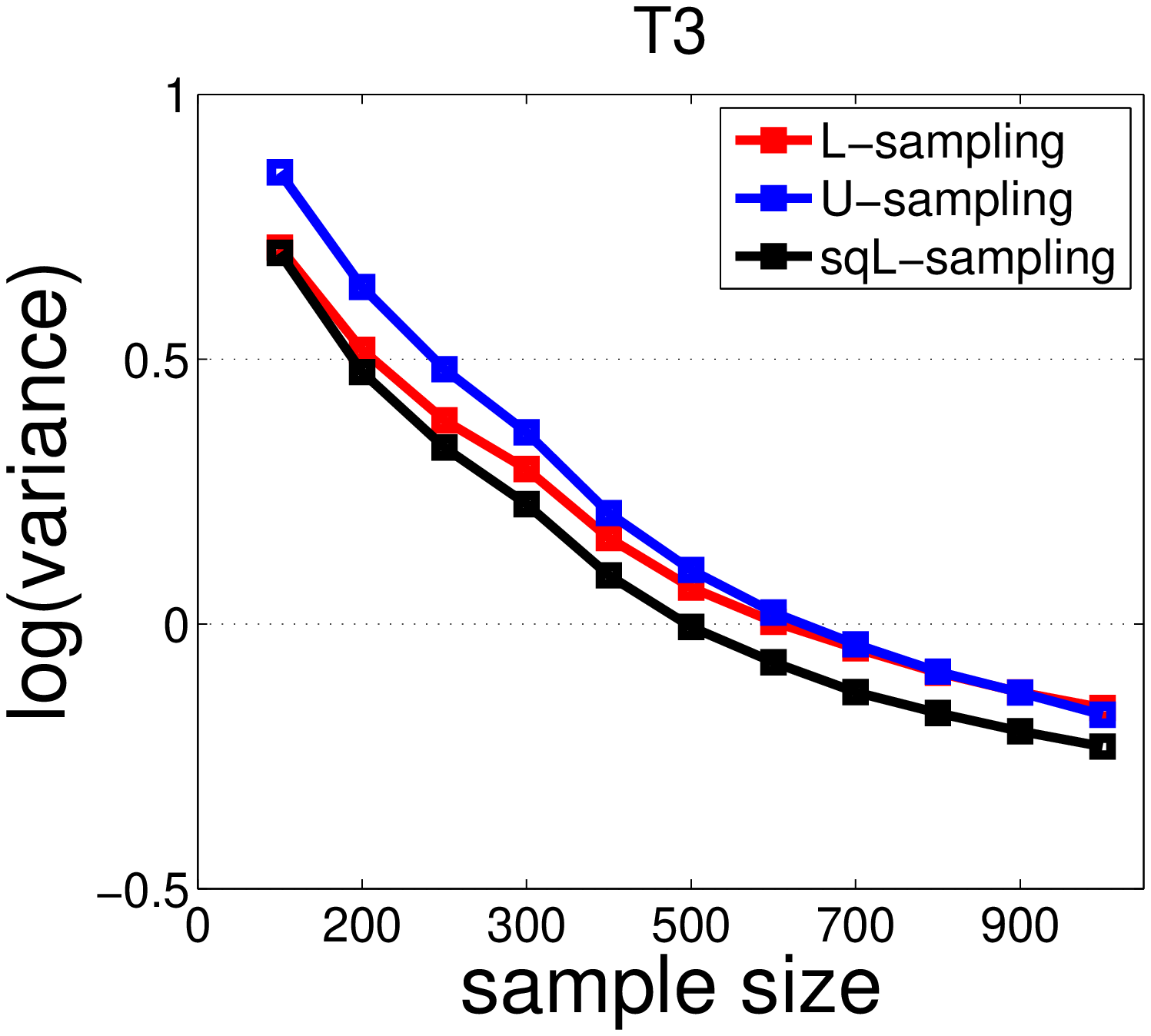}\hspace*{-0.1in}
\includegraphics[scale=0.17]{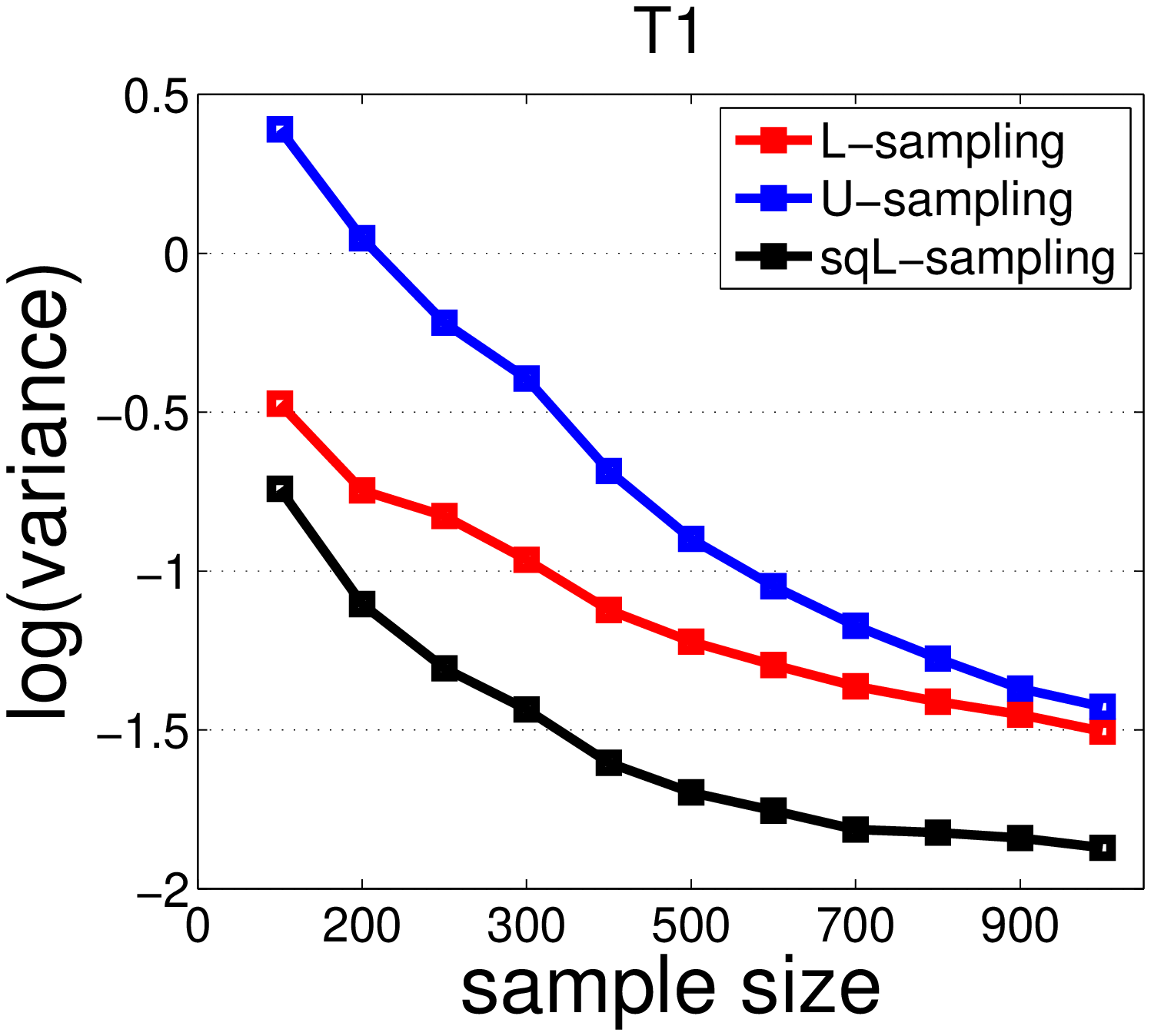}

\includegraphics[scale=0.17]{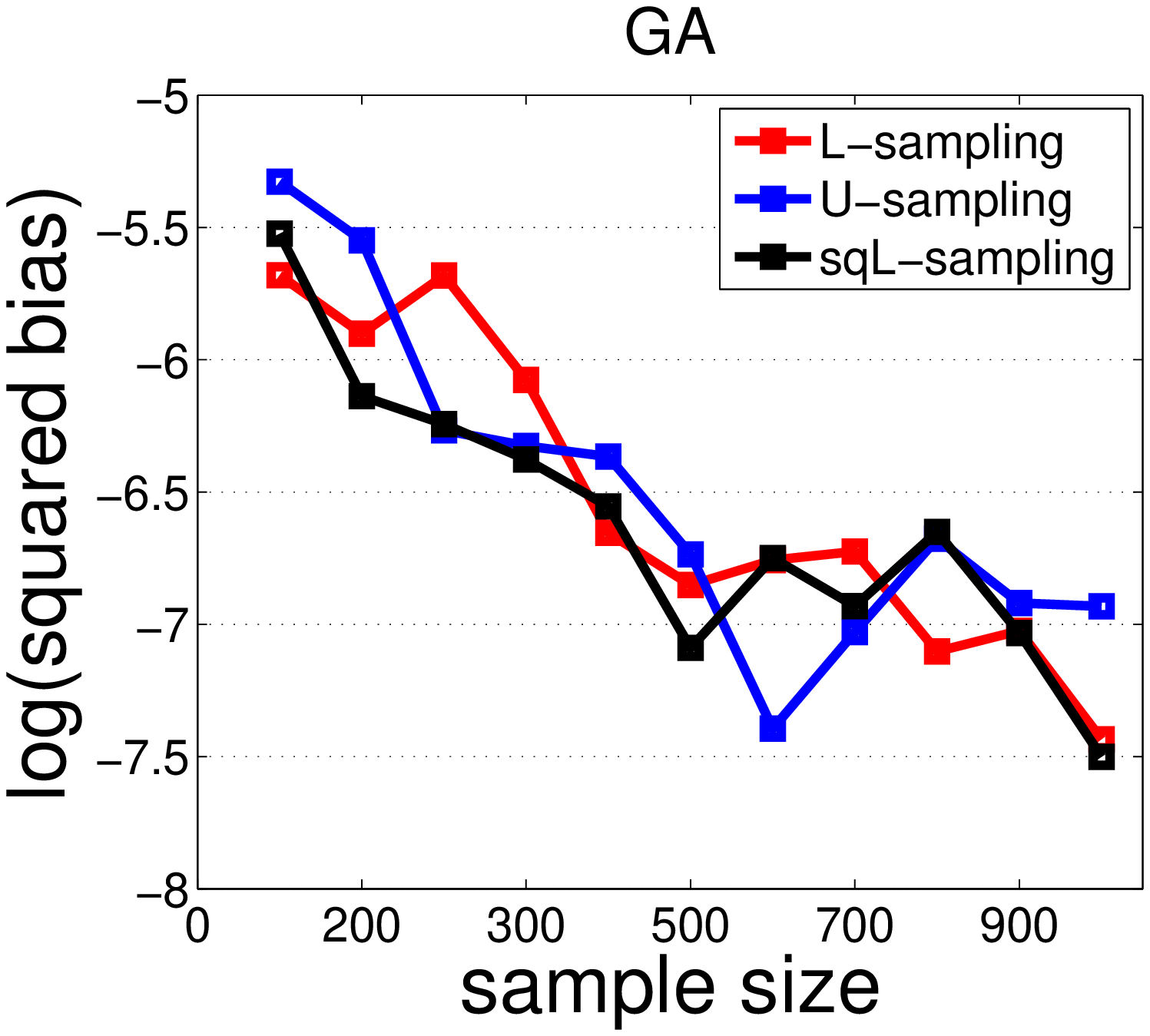}\hspace*{-0.1in}
\includegraphics[scale=0.17]{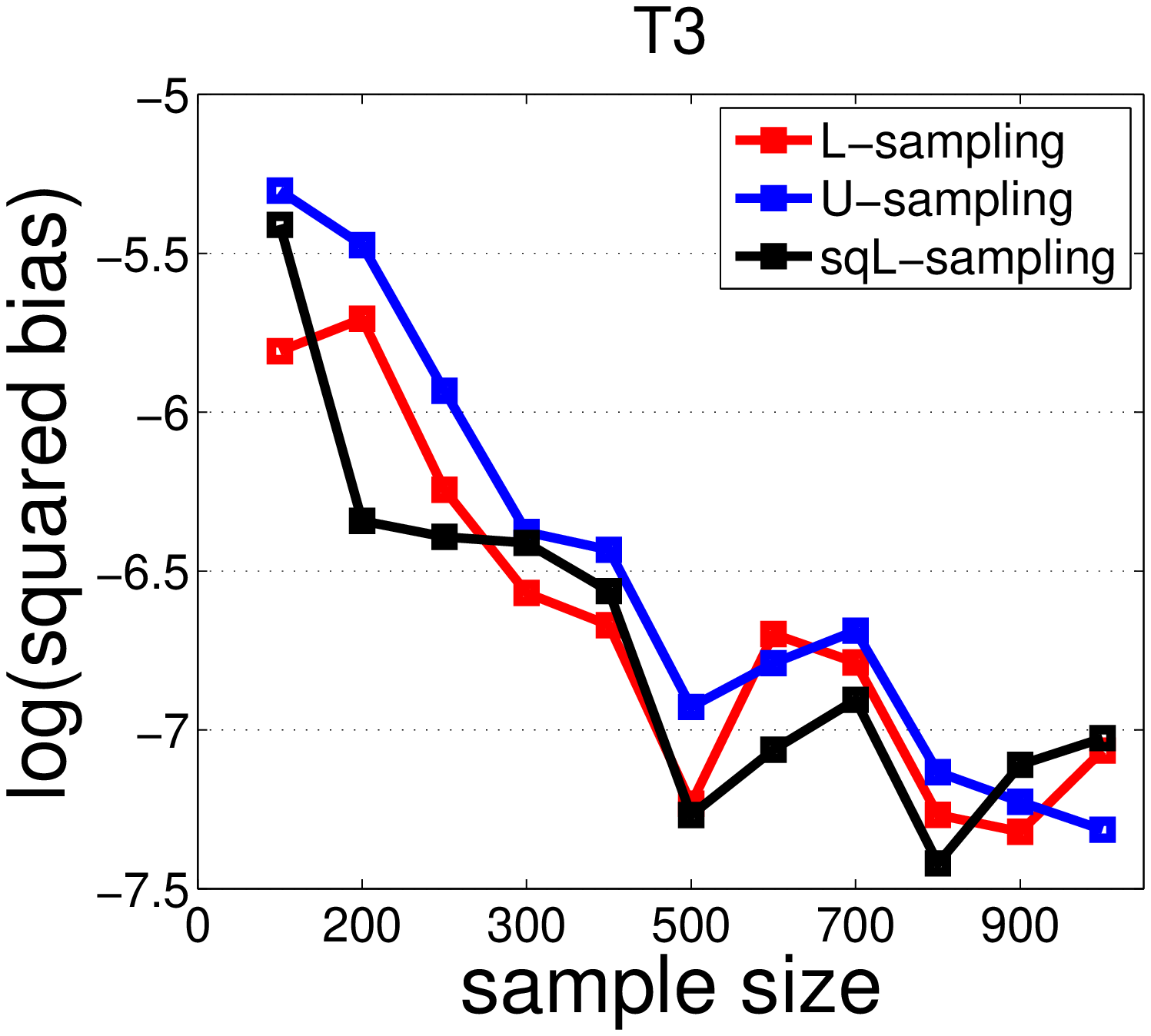}\hspace*{-0.1in}
\includegraphics[scale=0.17]{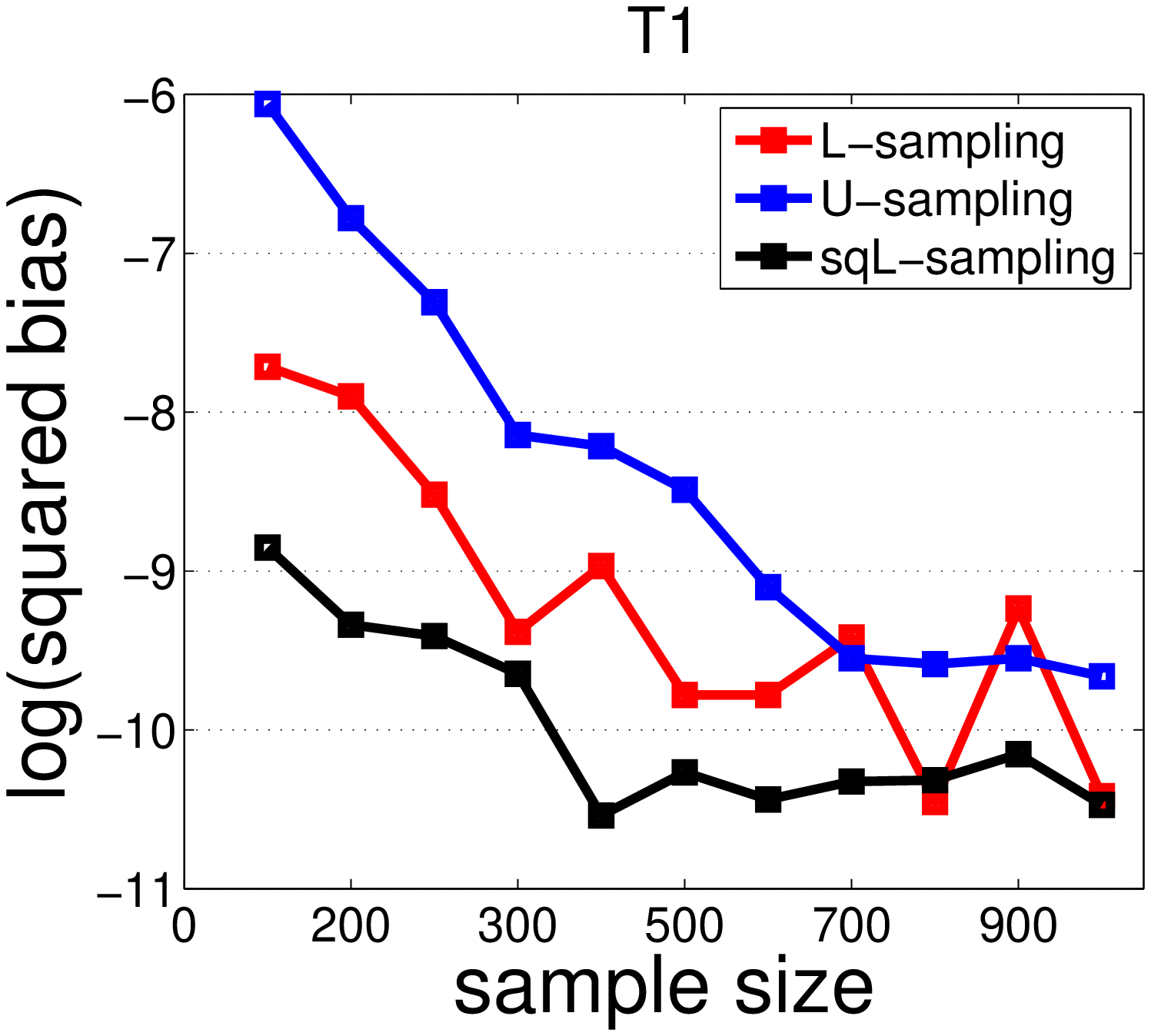}
\vspace*{-0.2in}
\caption{Comparison of variance and squared bias of the estimators  for different data, different samplings and different sample size.}\label{fig:3}
\vspace*{-0.1in}
\end{figure}
\begin{figure}[t]
\centering
\includegraphics[scale=0.17]{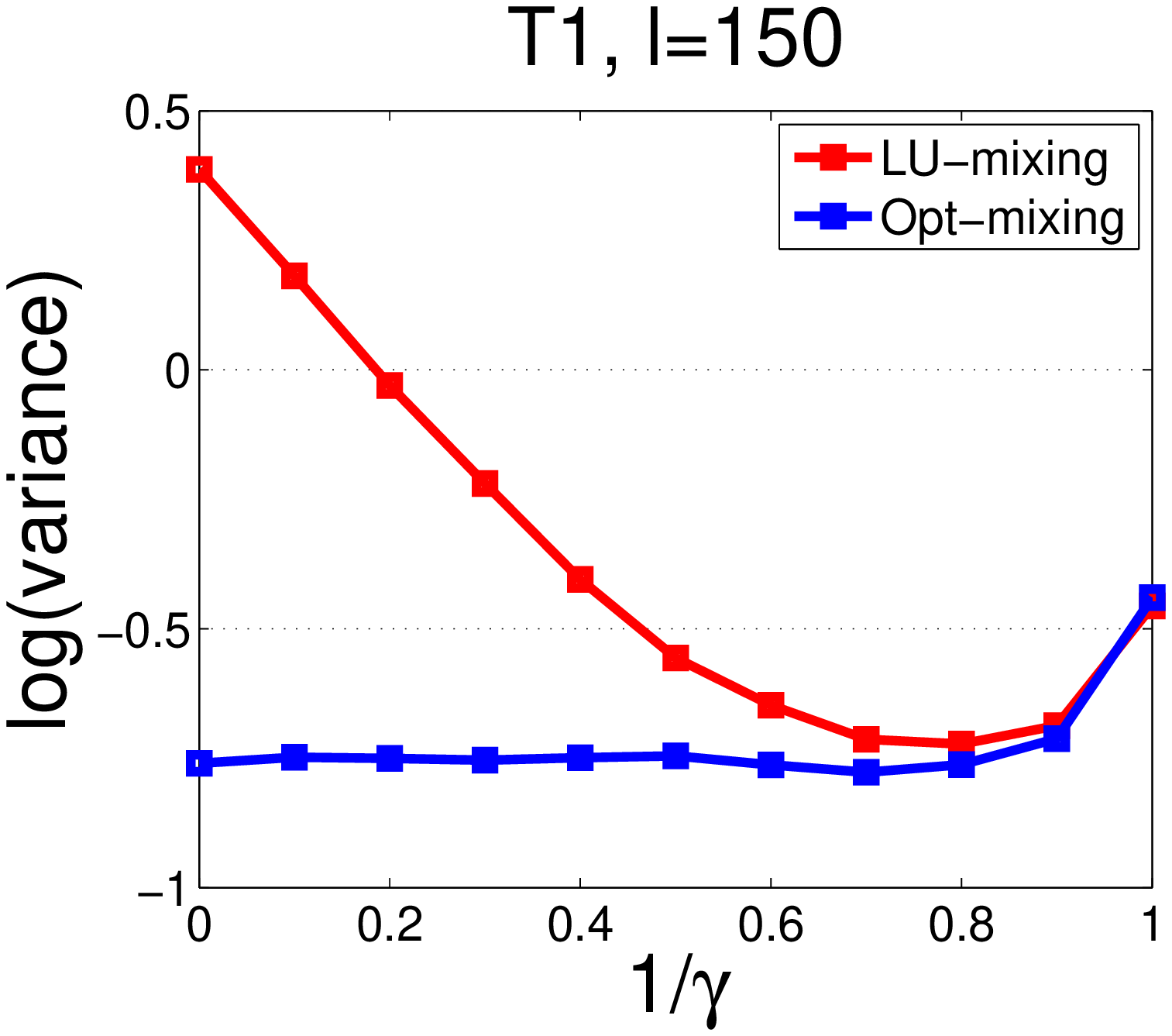}\hspace*{-0.1in}
\includegraphics[scale=0.17]{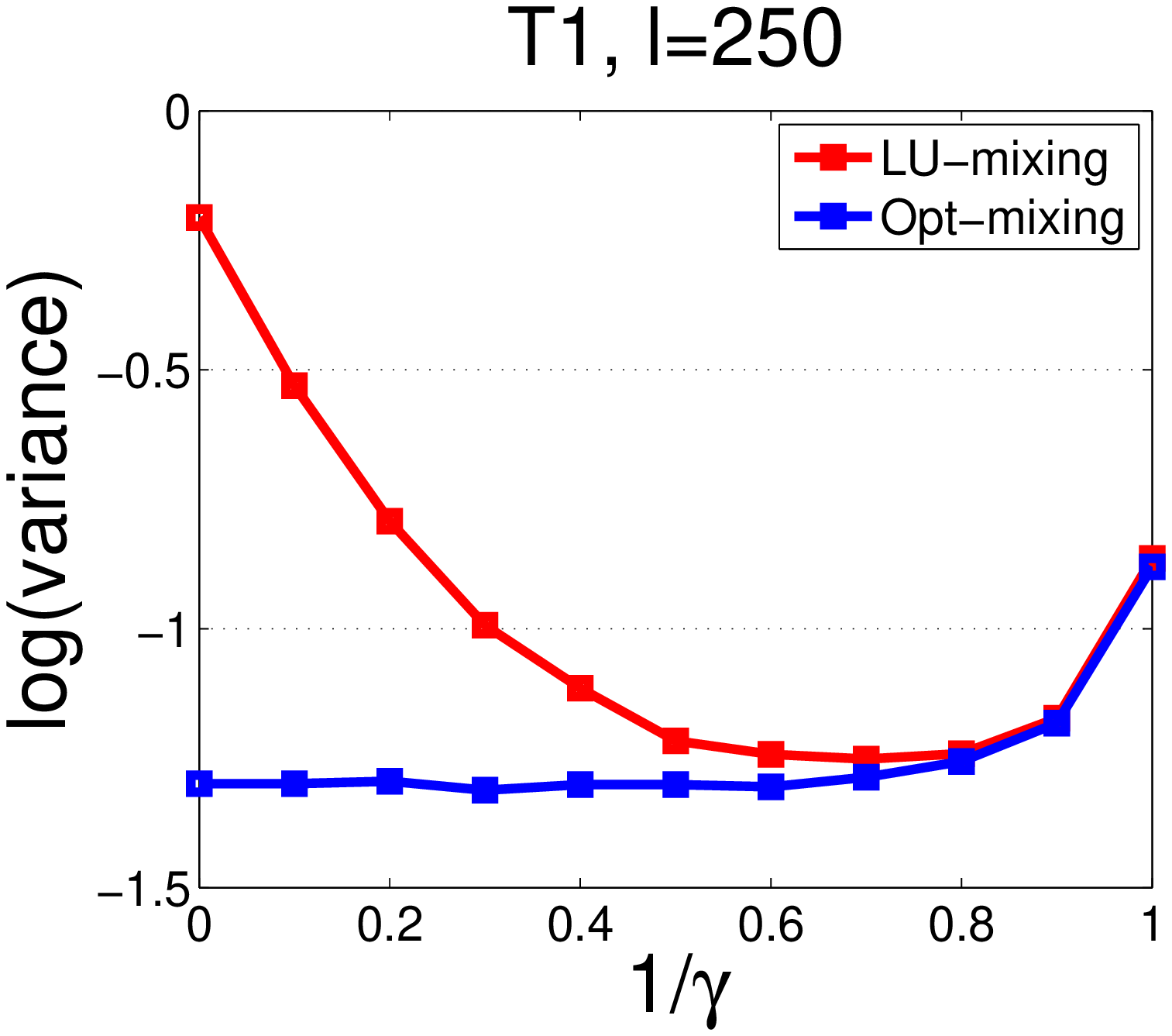}\hspace*{-0.1in}
\includegraphics[scale=0.17]{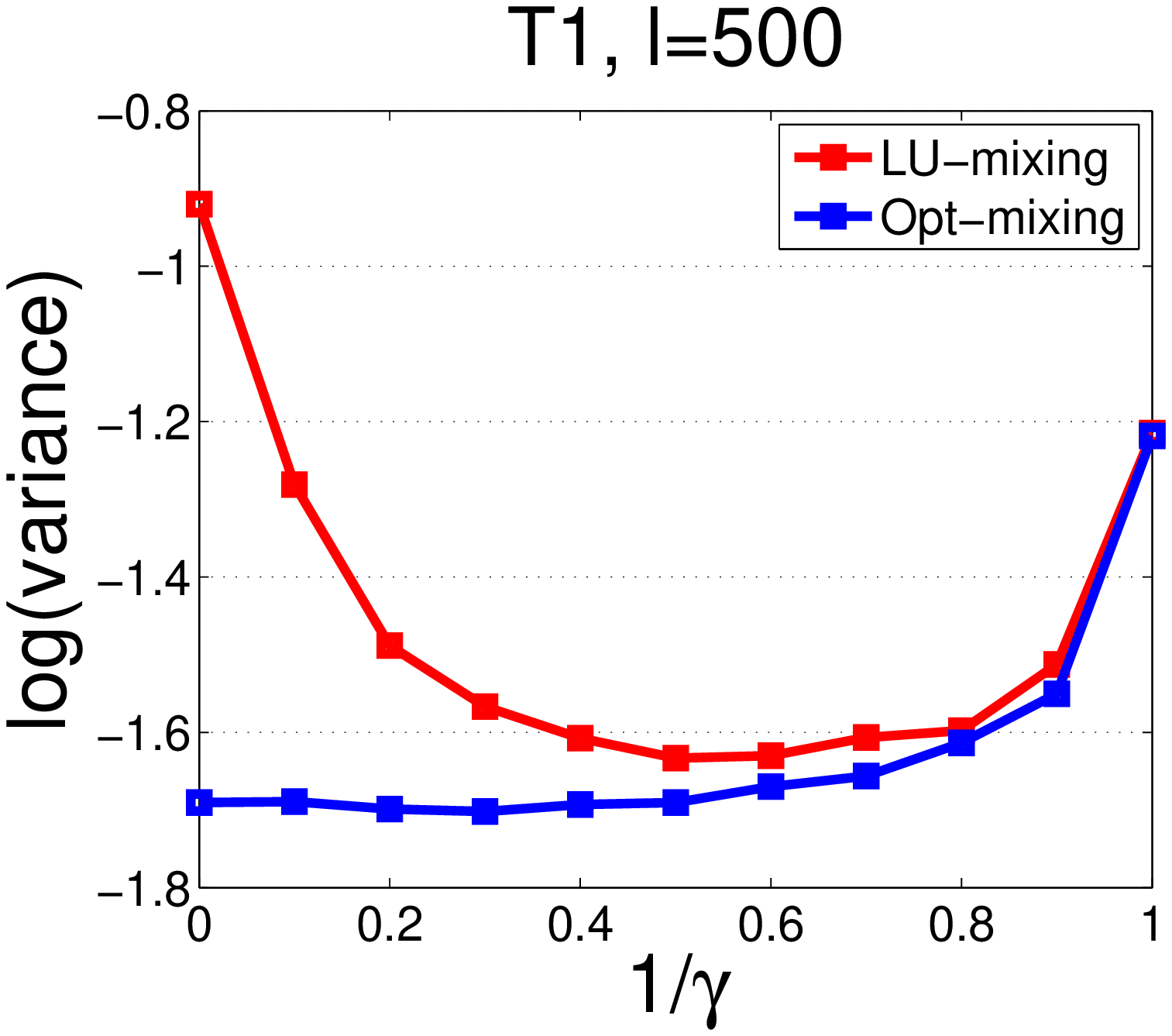}
\includegraphics[scale=0.17]{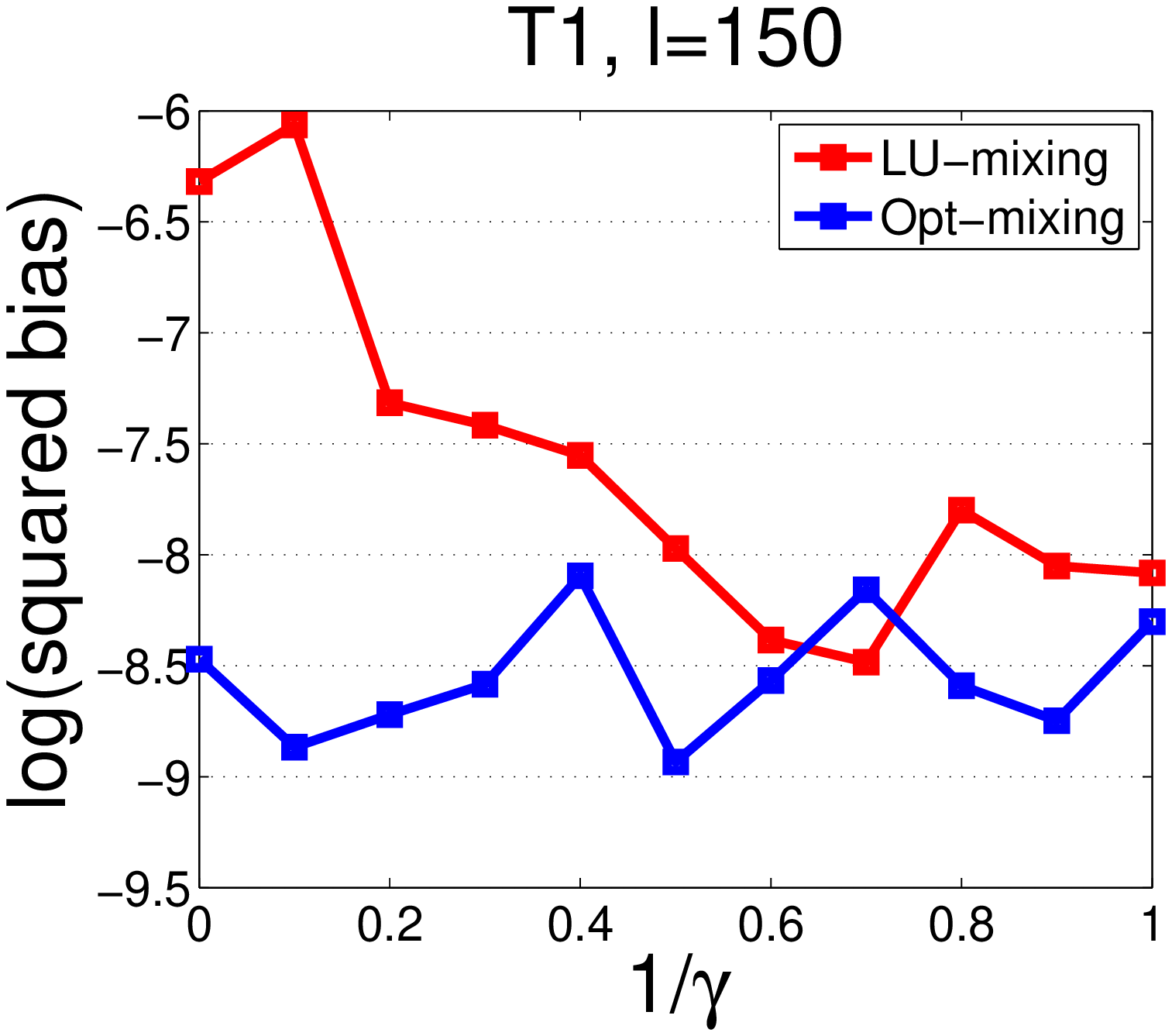}\hspace*{-0.1in}
\includegraphics[scale=0.17]{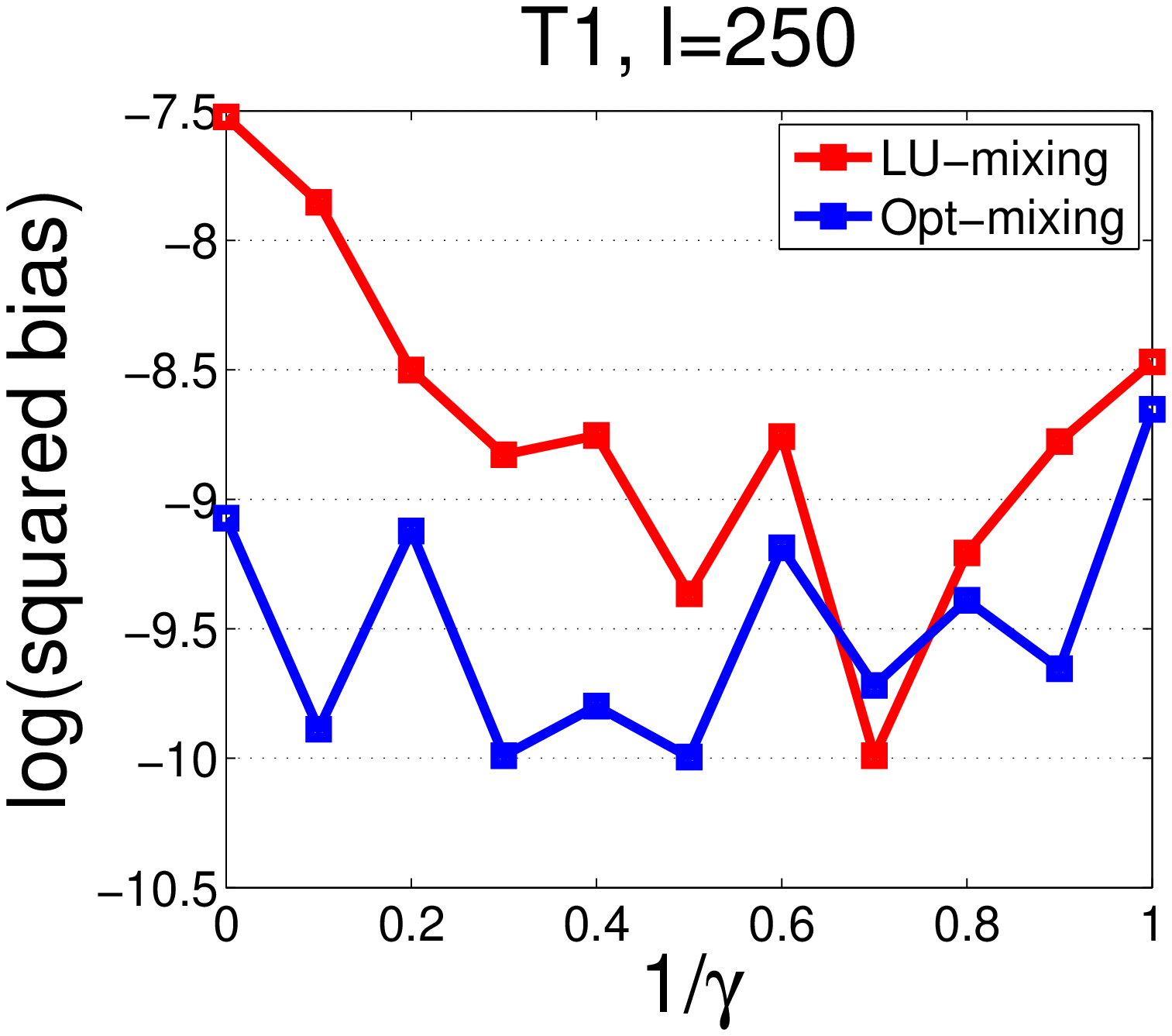}\hspace*{-0.1in}
\includegraphics[scale=0.17]{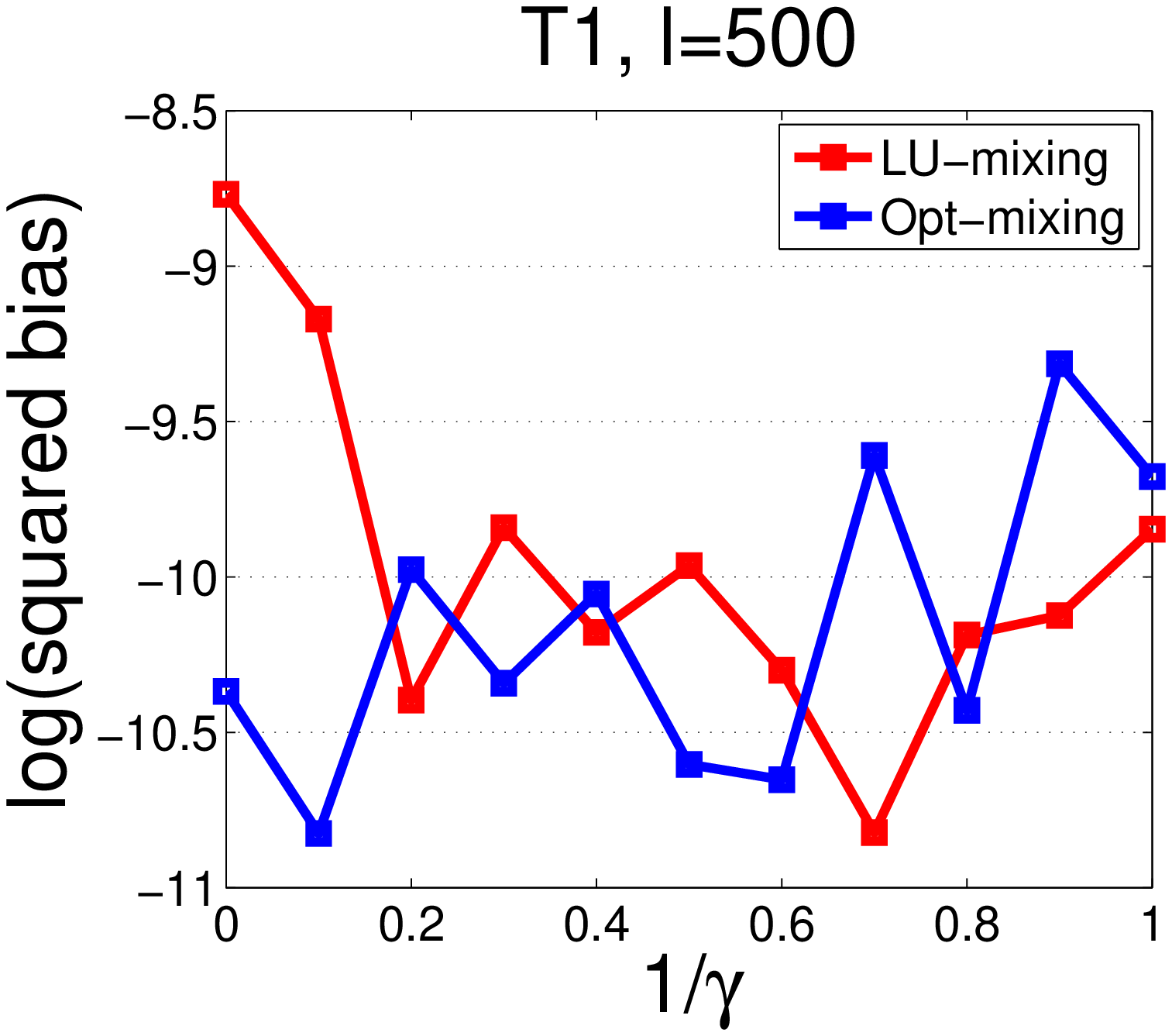}
\vspace*{-0.2in}
\caption{Comparison of variance and squared bias of the estimators  for different mixing strategies. Opt refers to our optimization based approach and LU refers to a convex combination of L-sampling and U-sampling with $\gamma^{-1}$ as the combination weight. }\label{fig:4}
\vspace*{-0.2in}
\end{figure}


Second, we compare the sampling scores found the constrained optimization with L-sampling and sqL-sampling. We vary the value of $\gamma$ from $1$ (corresponding to L-sampling) to $\infty$ (corresponding to sqL-sampling). A result with sampling size $\ell=500$ is shown in Figure~\ref{fig:2}. It demonstrate that intermediate samplings found by the proposed constrained optimization can perform better than both L-sampling and sqL-sampling. 

Finally, we apply CSS to over-constrained least square regression. To this end, we generate a synthetic data matrix $A\in\R^{m\times n}$ with $m=50$ and $n=1000$ similarly to~\cite{DBLP:conf/icml/MaMY14}. The output is generated by $y=A^{\top}\beta + \epsilon$ where $\epsilon\sim\mathcal (0, 9I_n)$ and $\beta=(\mathbf 1_{10}, 0.1\mathbf 1_{30}, \mathbf 1_{10})^{\top}$.  We compare the variance and bias of the obtained  estimators over $1000$ runs for different sampling distributions.  The results shown in Figure~\ref{fig:3} demonstrate the sqL-sampling gives  smaller variance and better bias of the estimators than L-sampling and U-sampling.  We also compare the proposed optimization approach with the simple mixing strategy~\cite{DBLP:conf/icml/MaMY14} that uses a convex combination of the L-sampling and the U-sampling. The results are shown in  Figure~\ref{fig:4}, which again support our approach.  

More results including relative error versus varying size $n$ of the target matrix,  performance on a real data set and the Frobenius norm reconstruction error  can be found in supplement. 

\section{Analysis}\label{sec:analysis}
In this section, we present major analysis of Theorem~\ref{thm:final-bound} and Theorem~\ref{thm:final-bound-2} with detailed proofs included in supplement. The key to our analysis is the following Theorem. 
\vspace*{-0.01in}
\begin{thm} \label{thm:bound-1} Let $Y, \Omega_1, \Omega_2$ be defined in ~(\ref{eqn:Y}). Assume that $\Omega_1$ has full row rank. We have
\[
\left\|A - P_YA\right\|_\xi^2 \leq \|\Sigma_2\|_\xi^2 + \left\|\Sigma_2\Omega_2\Omega_1^{\dagger}\right\|_\xi^2
\]
and 
\[
\left\|A - \Pi^2_{Y,k}(A)\right\|_\xi^2 \leq \|\Sigma_2\|_\xi^2 + \left\|\Sigma_2\Omega_2\Omega_1^{\dagger}\right\|_\xi^2
\]
where $\xi$ could be $2$ and $F$.
\end{thm}
\vspace*{-0.05in}
The first inequality was proved in~\cite{Halko:2011:FSR:2078879.2078881} (Theorem 9.1) and the second inequality is credited to~\cite{conf/focs/BoutsidisDM11} (Lemma 3.2)~\footnote{In fact, the first inequality is implied by the second inequality. }. 
Previous work on the spectral norm analysis also start from a similar inequality as above. They bound the second term by using $\|\Sigma_2\Omega_2\Omega_1^{\dagger}\|_2\leq \|\Sigma_2\Omega_2\|_2\|\Omega_1^{\dagger}\|_2$ and then bound the two terms separately.  However,  we will first write $\left\|\Sigma_2\Omega_2\Omega_1^{\dagger}\right\|_2=\|\Sigma_2\Omega_2\Omega_1^{\top}(\Omega_1\Omega_1^{\top})^{-1}\|_2$ using the fact $\Omega_1$ has full row rank, and then bound $\|(\Omega_1\Omega_1^{\top})^{-1}\|_2$ and $\|\Omega_2\Omega_1^{\top}\|_2$ separately. To this end, we will apply the Matrix Chernoff bound as stated in Theorem~\ref{thm:relative-bound} to bound $\|(\Omega_1\Omega_1^{\top})^{-1}\|_2$ and apply the matrix Bernstein inequality as stated in Theorem~\ref{thm:concentration} to bound  $\|\Omega_2\Omega_1^{\top}\|_2$. 

\begin{thm}[Matrix Chernoff~\citep{Tropp:2012:UTB:2347803.2347804}] \label{thm:relative-bound}
Let $\X$ be a finite set of PSD matrices with dimension $k$, and suppose that
 $ \max_{X \in \X} \lambda_{\max}(X) \leq B$.
Sample $\{X_1, \ldots, X_{\ell}\}$ independently  from $\X$. Compute
\begin{align*}
    \mu_{\max} = \ell \lambda_{\max}(\E[X_1]), \quad     \mu_{\min} = \ell \lambda_{\min}(\E[X_1])
\end{align*}
Then
\begin{align*}
&  \Pr\left\{\lambda_{\max}\hspace*{-0.05in}\left(\sum_{i=1}^{\ell} X_i\right) \geq (1 + \delta) \mu_{\max} \right\}\hspace*{-0.05in} \leq\hspace*{-0.05in} k \left[ \frac{e^{\delta}}{(1 + \delta)^{1 + \delta}}\right]^{\frac{\mu_{\max}}{B}} \\
& \Pr\left\{\lambda_{\min}\hspace*{-0.05in}\left(\sum_{i=1}^{\ell} X_i\right) \leq (1
  - \delta) \mu_{\min} \right\} \hspace*{-0.05in}\leq\hspace*{-0.05in} k \left[ \frac{e^{-\delta}}{(1 - \delta)^{1 - \delta}}\right]^{\frac{\mu_{\min}}{B}}
\end{align*}
\end{thm}
\vspace*{0.1in}
\begin{thm}[Noncommutative Bernstein Inequality~\citep{Recht:2011:SAM:1953048.2185803}] \label{thm:concentration}
Let $Z_1, \ldots, Z_L$ be independent zero-mean random matrices of dimension $d_1\times d_2$. Suppose $\tau_j^2 = \max\left\{\|\E[Z_jZ_j^{\top}]\|_2, \|\E[Z_j^{\top}Z_j\|_2\right\}$ and $\|Z_j\|_2 \leq M$ almost surely for all $k$. Then, for any $\epsilon > 0$,
\[
\Pr\left[\left\|\sum_{j=1}^L Z_j \right\|_2 > \epsilon \right] \leq (d_1+d_2)\exp\left[\frac{-\epsilon^2/2}{\sum_{j=1}^L \tau_j^2 + M\epsilon /3} \right]
\]
\end{thm}

Following  immediately from Theorem~\ref{thm:bound-1}, we have
\begin{align*}
\|A - P_YA\|_2& \leq \sigma_{k+1}\sqrt{1 + \|\Omega_2\Omega_1^{\top}(\Omega_1\Omega_1^{\top})^{-1}}\|^2_2\\
& \leq \sigma_{k+1}\sqrt{1 + \|\Omega_2\Omega_1^{\top}\|^2_2\lambda^{-2}_{\min}(\Omega_1\Omega_1^{\top})}\\
&\leq  \sigma_{k+1}(1 + \|\Omega_2\Omega_1^{\top}\|_2\lambda^{-1}_{\min}(\Omega_1\Omega_1^{\top})),
\end{align*}
where the last inequality uses the fact $\sqrt{a^2+b^2}\leq a+b$. Below we bound  $\lambda_{\min}(\Omega_1\Omega_1^{\top})$ from below and bound $\|\Omega_2\Omega_1^{\top}\|_2$ from above. 

\subsection{Bounding $\|(\Omega_1\Omega_1^{\top})^{-1}\|_2$} We will utilize Theorem~\ref{thm:relative-bound}  to bound $\lambda_{\min}(\Omega_1\Omega_1^{\top})$. Define $X_i = \v_i\v_i^{\top}/s_i$. It is easy to verify that
\begin{align*}
\Omega_1\Omega_1^{\top} = \sum_{j=1}^{\ell} \frac{1}{s_{i_j}}\v_{i_j} \v_{i_j}^{\top} = \sum_{j=1}^{\ell} X_{i_j}
\end{align*}
and $\E[X_{i_j}] = \frac{1}{\sum_{i=1}^n s_i}\sum_{i=1}^n s_iX_i = \frac{1}{k}I_k$, 
 where we use  $\sum_{j=1}^n s_j =k$ and $V_1^{\top}V_1=I_k$.
Therefore we have $\lambda_{\min}(\E[X_{i_j}]) = \frac{1}{k}$. Then the theorem below will follow Theorem~\ref{thm:relative-bound}. 

\begin{thm}\label{thm:omega1}
With a probability $1 - k\exp(-\delta^2\ell/[2kc(\s)])$, we have
\begin{align*}
\lambda_{\min}(\Omega_1\Omega_1^{\top}) \geq (1-\delta)\frac{\ell}{k}
\end{align*}
\end{thm}
Therefore, with a probability $1 - k\exp(-\delta^2\ell/[2kc(\s)])$ we have $\|(\Omega_1\Omega_1^{\top})^{-1}\|_2\leq \frac{1}{1-\delta}\frac{k}{\ell}$. 

\subsection{Bounding $\|\Omega_2\Omega_1^{\top}\|_2$} 
We will utilize Theorem~\ref{thm:concentration} to bound $\|\Omega_2\Omega_1^{\top}\|_2$. Define $Z_j = \u_{i_j}\v_{i_j}^{\top}/s_{i_j}$. Then $$\Omega_2\Omega_1^{\top} = \sum_{j=1}^{\ell} \frac{1}{s_{i_j}}\u_{i_j}\v_{i_j}^{\top}= \sum_{j=1}^lZ_j$$ and $\E[Z_j]=0$. In order to use the matrix Bernstein inequality, we will  bound $\max_i\|Z_i\|_2 =\max_i \frac{\|\u_i\v_i^{\top}\|_2}{s_i} \leq q(\s)$ and $\tau_j^2\leq \frac{(\rho +1 - k)c(\s)}{k}$. Then we can prove the following theorem. 
\begin{thm}\label{thm:cross}
With a probability $1 - \delta$, we have
\[
\|\Omega_2\Omega_1^{\top}\|_2\leq \sqrt{2c(\s)\frac{(\rho+1-k)\ell\log(\frac{\rho}{k})}{k}} + \frac{2q(\s)\log(\frac{\rho}{k})}{3}.
\]
\end{thm}
We can complete the proof of Theorem~\ref{thm:final-bound} by combining the bounds for $\|\Omega_2\Omega_1^{\top}\|_2$  and $\lambda^{-1}_{\min}(\Omega_1\Omega_1^{\top})$ and by setting  $\delta=1/2$ in Theorem~\ref{thm:omega1} and using union bounds.


\section{Discussions and Open Problems}
From the analysis, it is clear that the matrix Bernstein inequality is the key to derive the sampling dependent bound for $\|\Omega_2\Omega_1^{\top}\|_2$. For bounding $\lambda_{\min}(\Omega_1\Omega_1^{\top})$, similar analysis using matrix Chernoff bound has been exploited before for randomized matrix approximation~\cite{gittens-2011-spectral}.  

Since Theorem~\ref{thm:bound-1} also holds for the Frobenius norm, it might be interested to see whether we can derive a sampling dependent Frobenius norm error bound that depends on $c(\s)$ and $q(\s)$, which,  however, still  remains as an open problem for us. Nonetheless, in experiments (included in the supplement) we observe similar phenomena about the performance of L-sampling, U-sampling and sqL-sampling.    

Finally, we briefly comment on the analysis for least square approximation using CSS. Previous results~\cite{journals/siammax/DrineasMM08,conf/soda/DrineasMM06,Drineas:2011:FLS:1936922.1936925} were built on the structural conditions that are characterized by two inequalities 
\begin{align*}
&\lambda_{\min}(\Omega UU^{\top}\Omega)\geq 1/\sqrt{2}\\
& \|U^{\top}\Omega^{\top}\Omega U^{\perp}{U^{\perp}}^{\top}\b\|^2_2\leq\frac{\epsilon}{2}\|U^{\perp}{U^{\perp}}^{\top}\b\|^2_2
\end{align*}
The first condition can be guaranteed by Theorem~\ref{thm:omega1} with a high probability. For the second condition, if we adopt a worse case analysis
\begin{align*}
 \|U^{\top}\Omega^{\top}\Omega U^{\perp}{U^{\perp}}^{\top}\b\|^2_2\leq  \|U^{\top}\Omega^{\top}\Omega U^{\perp}\|^2_2\|{U^{\perp}}^{\top}\b\|^2_2
\end{align*}
and bound the first term in R.H.S of the above inequality using Theorem~\ref{thm:cross}, we would end up with a worse bound than existing ones that bound the left term as a whole. Therefore the naive combination can't yield a good sampling dependent error bound for the approximation error of least square regression. 

\vspace*{-0.1in}
\section{Conclusions}
In this paper, we have presented a sampling dependent spectral error bound for CSS. The error bound brings a new distribution with sampling  probabilities proportional to the square root of the statistical leverage scores and exhibits more tradeoffs and insights than existing error bounds for CSS.  We also develop a constrained  optimization algorithm with an efficient bisection search to find better sampling probabilities for the  spectral norm reconstruction. Numerical simulations demonstrate that the new sampling distributions lead to improved performance. 

\bibliography{all}
\bibliographystyle{icml2015}

\end{document}